\documentclass[12pt]{article}
\usepackage{amsfonts}
\usepackage{amsmath, amsthm, amssymb}
\newcommand {\eG}  %usato
        { {\cal{E}}^{'} } 
\newcommand {\betapi} %usato
        {\beta^{\pi}}

\newcommand {\scalaires} %usato
        {\mathbb{K}}

\usepackage[latin1]{inputenc}
\usepackage[T1]{fontenc}
\DeclareMathAlphabet\EuFrak{U}{euf}{m}{n}%BOLD EULER FRAKTUR
\SetMathAlphabet\EuFrak{bold}{U}{euf}{b}{n} % GOTIC FONT
\newtheorem{tm}{Theorem}[section]

\newtheorem{Def}{Definition.}[section]
\newtheorem{Cor}{Corollary}[section]
\newtheorem{Oss}{Remark}[section]
\newtheorem{ex}{Example}[section]
\newtheorem{lemma}{Lemma}[section]
\newtheorem{convenzione}{Convention}[section]
\newtheorem{notazione}{Notation}[section]
\newtheorem{assume}{Hypothese}[section]
\newcommand {\dimostrazione}
        {{\sc Proof}.}
\newcommand {\halmos}
        {\ensuremath{\raisebox{-.3ex}{\rule{1ex}{1em}}}}
\newcommand {\gG} % g gotique}
        {\EuFrak {g} }
\newcommand {\hG} % h gotique}
        {\EuFrak {h} }
\newcommand {\qG} % h gotique}
        {\EuFrak {q} }
\renewcommand{\L}
        {{\partial}}
\newcommand {\ad}
        {{\rm ad}}
\newcommand {\Mult}
        {\rm{Mult}}
\newcommand {\campo}
        {\mathbb{K}}

\newcommand {\Z}
        {\mathbb{Z}}

\newcommand {\Naturali} %
        {\mathbb{N}}  
\newcommand {\Razionali} %
        {\mathbb{Q}}  
\title{Universal representations of Lie algebras by coderivations}
\author{ 
Emanuela Petracci\footnote{Reserch partially supported by the french association 
A.F.F.D.U.}
        }
\date{}
\begin{document}
\maketitle
\abstract{
A class of representations of  a Lie superalgebra (over a commutative 
superring) in its symmetric algebra is studied. 
As an application we get  a direct and 
natural proof of a strong form of the Poincar{\'e}-Birkhoff-Witt theorem,  
extending this 
theorem to a class of nilpotent Lie superalgebras.  Other applications 
are 
presented. Our results are new already for Lie algebras.}

%\tableofcontents
\section{Introduction}
We consider a commutative superring $\campo= \campo_0 +\campo_1$,
and for 
all Lie superalgebra $\gG$ over $\campo$ we
study the representations of $\gG$ on its symmetric algebra $S(\gG)$ 
which are {\it by coderivations and universal}. 

The symmetric  algebra $S(\gG)$ has a natural structure 
of coalgebra, so
we have a notion of coderivation of $S(\gG)$. A representation $\rho$ of
$\gG$ in $S(\gG)$ 
is called {\it by coderivations} if $\rho(a)$ is a coderivation of 
$S(\gG)$ for all $a\in\gG$. 
We focus on representations $\rho$ by coderivations which are 
{\it universal}. 
This means informally that $\rho$ is given by 
a formula independent of $\gG$ (see definition \ref{def32}). 

To each  formal power series 
$
\varphi = c_0+c_1t+\cdot\cdot\cdot\in\campo_0[[t]]
$ 
we associate   a family 
$\Phi(a)\equiv\Phi^a$ of
coderivations of $S(\gG)$ depending linearly of $a\in\gG$ (see formula 
(\ref{terzaL})). 

We show that  $\Phi$ is an universal representation by
coderivations if and only if
\begin{equation}
\label{prima}
  \varphi(t)
  \frac{\varphi(t+u) - \varphi(u)}{t}
+ \varphi(u)
  \frac{\varphi(t+u)-\varphi(t)}{u}
+ \varphi(t+u)
= 0
\end{equation}
in $\campo_0[[t,u]]$. 

We show that, for Lie algebras over a $\Razionali$-algebra, all 
 universal representations are of this form (theorem 
\ref{tm72}). 

The most interesting case is when the constant term $c_0$ is equal to $1$. 
In this case, it is a simple matter to solve the functional equation,
because we show (theorem \ref{lemma23}) that it is equivalent to the 
functional equation for the
exponential function. This last equation has a non-trivial solution
exactly when $\campo$ contains $\Razionali$. In this case the unique
solution of (\ref{prima}) is the generation function for Bernoulli numbers 
\begin{eqnarray*}
\label{seconda}
       \varphi
     = \varphi_1
\equiv \frac{t}{e^t-1}
.
\end{eqnarray*}

Let $N\geq 2$ be an integer. If we restrict to an $N$-nilpotent Lie
superalgebras $\gG$ over $\campo$, we get similar results. To a truncated 
power series $\varphi\in\campo_0[[t]]/t^{N}$ is associated a family of 
coderivations $\Phi^a$ depending linearly of $a\in\gG$. We  show  that $\Phi$ 
is an universal  representation by coderivations if and only if $\varphi$ 
verifies equation (\ref{prima}) in $\campo_0[[t,u]]/I_N$, where $I_N$ is the 
ideal generated by $\{ t^iu^{N-1-i}| 0\leq i\leq N-1 \}$. 
There exists a solution with $ \varphi(0)=1 $ exactly when 
$
  \frac{1}{2}
 , \frac{1}{3}
 , ...
 , \frac{1}{N}
\in\campo
$,  
and in this case the unique solution is $\varphi_1 \mod t^N$.
\newline

We explain  the relation of these results with the 
Poincar{\'e}-Birkhoff-Witt theorem. Let $U(\gG)$ be the enveloping algebra of
$\gG$ and assume that $\campo\supseteq\Razionali$. We use the 
representation obtained using
the function $\frac{t}{e^t-1}$ to define a symbol map
$\sigma: U(\gG)\rightarrow S(\gG)$. We show that $\sigma$ is an 
inverse
for the symmetrization $\beta: S(\gG)\rightarrow U(\gG)$, which gives  a 
natural and direct proof of the fact that  $\beta$ is an 
isomorphism. 

Let  $\gG$ be a $N$-nilpotent Lie superalgebra over a commutative 
superring containing    $ \frac{1}{2}, \frac{1}{3}, ..., \frac{1}{N}$. 
Also in this case there is  a canonical symbol map $\sigma: 
U(\gG)\rightarrow S(\gG)$ which is an isomorphism. For $N=2$ this is due to M. El-Agawany and  A. 
Micali ({\bf\cite{ElM}}). The case $N\geq 3$ is new.
\newline

The equation (\ref{prima}) is a particular case of an equation studied in 
section \ref{sec22}. As an application of this "more general equation" we 
study  the universal
representations by coderivations of  $\gG\times\gG$ on $S(\gG)$. 
\newline

Let $\gG$ be any Lie superalgebra over a superring. The enveloping 
algebra $U(\gG)$ also has a natural structure of 
coalgebra. Assume that the Poincar{\'e}-Birkhoff-Witt theorem is 
verified. Using the isomorphism $\sigma^{-1}$, we show 
that an universal representation by coderivations  
$\Phi: \gG\rightarrow End(S(\gG))$ gives an universal representation by 
coderivations $F: \gG\rightarrow End(U(\gG))$. We get a family of 
representations interpolating the left and the right regular 
representation and the adjoint representation of $\gG$ in $U(\gG)$.
\newline
\newline
{\it Acknowledgements.}
The present paper is extracted from my PhD thesis ({\bf\cite{Pet}}). I 
would like  to thank Michel Duflo for having been my advisor.

\section{Lie superalgebras over a superring}
In this section we recall the basic definitions and examples used in the
text, they are  from super linear algebra ({\bf\cite{Lei}}).

We say that $\campo$ is a {\it superring} if it is an unitary ring
graded over $\Z/2\Z$. We denote by $\campo_0$ 
and $\campo_1$ the subgroups of elements 
with even and odd degree,  for each not-zero homogeneous
element $a\in \campo$ we denote by  
$p(a)$ its degree. We have $1\in \campo_0$. 

The superring $\campo$ is called {\it commutative} if 
$a b = (-1)^{p(a)p(b)}b a$ for all
homogeneous and not-zero $a,b\in\campo$, and 
$a^2 = 0$ for $a\in \campo_1$.

\begin{convenzione}
Each time we use the symbol $p(a)$ for an element $a$ of a graded 
group 
occurring in a linear expression, it is implicitly assumed that it is not-zero 
and homogeneous. Moreover the expression is extended by linearity. 
For example, the expression above will be written as 
$ab=(-1)^{p(b)p(a)}ba$ for any $a,b\in \campo$.

\end{convenzione}

From now to the end of this section, $\campo$ will be a fixed commutative 
superring. 
\newline
We denote by $\campo^\times_0\subseteq\campo$  the subgroup of 
 invertible elements of $\campo_0$.

\begin{Def}
A commutative  group 
$(M,+)$ graded 
over 
$\Z/2\Z$ is a 
$\campo$-module if it is equipped with a bilinear  application 
$M\times\campo\rightarrow M$ such that, for any $\alpha,\beta\in 
\campo$ 
and $m,n\in M$ we have 
\begin{eqnarray*}
(m\alpha )\beta = m(\alpha\beta )
,\ \ \ 
p(m\alpha) = p(m) + p(\alpha)
.
\end{eqnarray*}
We denote by $M_0$ and $M_1$ the $\campo_0$-submodules composed  of even 
and odd elements. 
\end{Def}
\noindent
In a $\campo$-module $M$ we use the notation  $\alpha 
m:=m\alpha(-1)^{p(\alpha)p(m)}$, for any $m\in M$ and $\alpha\in 
\campo$. Let $N$ be another $\campo$-module. A  map $f:M\rightarrow N$ is 
a {\it morphism of $\campo$-modules} 
if $f(m\alpha)=f(m)\alpha$ for any $m\in M$ and $\alpha\in\campo$.

\begin{Def}
We say that $A$ is a 
$\campo$-superalgebra if it a $\campo$-module equipped of a
distributive application $A\times A\stackrel{\cdot}{\rightarrow} A$ such that 
\begin{eqnarray*}
p(a\cdot b) = p(a)+p(b)
, \ \ \ \ 
(a\cdot b)\alpha 
= a\cdot(b\alpha)
= (-1)^{p(b)p(\alpha)}
  (a\alpha)\cdot b
\end{eqnarray*}
for any $a,b\in A$ and $\alpha\in\campo$. 
We say that $A$ is commutative if $a\cdot b = (-1)^{p(a)p(b)}b\cdot a$ 
for $a,b\in A$, and $c^2 = 0$ for $c\in A_1$.
\end{Def}
\noindent
Let $A$ and $B$ be two $\campo$-superalgebras. A map 
$f: A\rightarrow B$ is said a {\it morphism of $\campo$-superalgebras} if 
it is a morphism of $\campo$-modules such that 
\begin{eqnarray*}
  p(f(a)) =p(a)
, \ \ \ 
  f(a\cdot b)=f(a)\cdot f(b)
, \ \ \ \ \ 
  \forall a,b\in A
.
\end{eqnarray*}

\begin{notazione}
\label{not21}
Let $A$ be a $\campo$-superalgebra and $a\in A$. We denote by $a^L : 
A\rightarrow A$ the
left multiplication by $a$, and by $a^R : A\rightarrow A$ the right
multiplication by $a$: $a^R(b) = (-1)^{p(a)p(b)}b\cdot a$, for any $b\in
A$.
\end{notazione}
\noindent
The following is  our definition of Lie superalgebra.

\begin{Def}
\label{Defin13}
Let   $\gG$ be a $\campo$-superalgebra such 
that its product $[\cdot,\cdot ] : \gG\times\gG\rightarrow\gG$ 
verifies
\begin{eqnarray}
\label{labi}
&& \left[X,Y\right]
= -(-1)^{p(X)p(Y)}
  [Y,X]
, \ \ \forall X,Y\in\gG
\\
&&\label{labiii}
[X,X]=0
, \ \ \forall X\in\gG_0
\\
\label{labii}
&&  [[X,Y],Z]
= [X, [Y,Z]] 
- (-1)^{p(Y)p(X)}
  [Y, [X,Z]]
,\ \ \forall X,Y,Z\in\gG
\\
&&\label{labiv}
[Y, [Y,Y]]=0
, \ \ \forall Y\in\gG_1
.
\end{eqnarray}
Such $\gG$ is called a Lie $\campo$-superalgebra.
\end{Def}
\noindent
The product in a Lie superalgebra is called Lie product or Lie
bracket, and (\ref{labii}) is  the {\it Jacobi identity}.
\begin{Oss}
If $2\in\campo$ is invertible (\ref{labiii}) follows from (\ref{labi}).
If $3\in\campo$ is invertible (\ref{labiv}) follows from (\ref{labi}) and 
(\ref{labii}).
\end{Oss}
\noindent
As explained in 
{\bf \cite{BMP}}, if $\gG_1\neq \{0\}$ and $2\in\campo$ is not  
invertible, the definition \ref{Defin13} is not  the right one, but it is 
sufficient for the porpoise of this text.

We end this section with some useful examples.

\begin{ex}
If $A$ is a commutative superring, $A[[z]]$ denotes the set formal
series in
$s$, with coefficients in $A$. It inherits the graduation
$(A[[z]])_0 = A_0[[z]]$, $(A[[z]])_1 = A_1[[z]]$ and a natural structure
of commutative
superring.
\end{ex}

\begin{ex}
Let $M,N$ be two $\campo$-modules.

a) $Hom(M,N)$ is  the group of  functions $F:M\rightarrow N$
which are morphisms of $\campo$-modules. It is graded in the following 
way: $F$ is even if $F(M_0)\subseteq N_0$ and $F(M_1)\subseteq N_1$, $F$ 
is odd if $F(M_0)\subseteq N_1$ and $F(M_1)\subseteq N_0$. 
Moreover, $Hom(M, N)$ is a $\campo$-module by 
$ F\alpha: v
  \mapsto (-1)^{p(v)p(\alpha)}F(v)\alpha
$, 
  for all  $\alpha\in\campo$, $v\in M$.

b) $M\otimes N$ is the $\campo$-module generated by 
$\{v\otimes w; v\in M,w\in N, \}$
with relations
\begin{eqnarray*}
  &&(v_1+v_2)\otimes w
= v_1\otimes w + v_2\otimes w
 \\  
  &&v \otimes (w_1+w_2)
= v\otimes w_1+ v\otimes w_2
\\
&&(v\otimes w)\alpha
= v\otimes w\alpha
= (-1)^{p(w)p(\alpha)}
  v\alpha\otimes w 
,\ \
\forall\alpha\in\campo
\end{eqnarray*}
and graduation  $p(v\otimes w)= p(v)+p(w)$.

c) The tensor algebra of $M$ is 
$ T(M) 
:= \campo
+ (M\otimes M)
+ (M\otimes M\otimes M)
+\cdot\cdot\cdot
$ 
with product 
$     (v_1\otimes \cdot\cdot\cdot \otimes v_i)
\cdot (v_{i+1}\otimes \cdot\cdot\cdot \otimes v_n)
      = v_1\otimes \cdot\cdot\cdot \otimes v_n
$, for all $i,n\in\Naturali$. It is an associative  $\campo$-superalgebra.
\end{ex}

\begin{ex}
Let  $M$ a $\campo$-module, 
$Hom(M,M)$ is a Lie $\campo$-superalgebra with bracket 
$[F,G]=F\circ G - (-1)^{p(F)p(G)} G\circ F$, for any $F,G\in Hom(M,M)$.
\end{ex}

\section{Symmetric algebras }
\label{sez12}
In all this section $\campo$ is a commutative superring. 

Let  $M$  be a $\campo$-module, we recall the definition of its symmetric
algebra $S(M)$. 
The tensor algebra $T(M)$ 
 contains the ideal $I$ generated by 
$$ \left\{
    v\otimes w - (-1)^{p(v)p(w)}w\otimes v
  , u\otimes u
  | \ v,w\in M, u\in M_1
  \right\}
,
$$ 
and we  define $S(M):= T(M)/I$. It is a commutative and associative 
$\campo$-superalgebra. 
We have 
$ S(M)
= \campo\oplus \bigoplus_{n=1}^{\infty}
  S^n(M)
$, 
where  $S^n(M)$ is the $\campo$-module generated by products of $n$ 
 elements of $M$. 

\subsection{Formal functions}
\label{sec121}
We recall that $S(M)$ has a natural structure of
cocommutative Hopf superalgebra, and in particular it is a coalgebra. This 
means that $S(M)$ is
equipped of three  morphisms of superalgebras 
$\Delta : S(M)\rightarrow S(M)\otimes S(M)$, 
$\epsilon : S(M)\rightarrow \campo$, 
$\delta: S(M)\rightarrow S(M)$, 
such that
\begin{eqnarray} 
\label{lzero}
&&      (id\otimes\Delta)
\circ \Delta
    = (\Delta\otimes id)
\circ \Delta
\\
\label{assiomaq}
&&{\Mult}\circ (id\otimes \delta)\circ \Delta
= {\Mult}\circ(\delta\otimes id)\circ\Delta
= \epsilon
\\
&&
\label{assiomaL}
      {\Mult} 
\circ (id\otimes \epsilon)\circ \Delta 
    = {\Mult}
\circ (\epsilon\otimes id)\circ \Delta 
    = id
\\
&& 
\label{assiomap}
  \Delta
= \sigma\circ\Delta
\end{eqnarray} 
where ${\Mult} :S(M)\otimes S(M){\"U}\rightarrow S(M)$ is
the multiplication of $S(M)$,  
$\sigma : S(M)\otimes S(M)\ni W\otimes Z  \mapsto (-1)^{p(W)p(Z)} Z\otimes W
$ 
is the exchange operator. We call $\delta$ an antipode, and each even 
morphism of $\campo$-modules verifying (\ref{lzero}) is called an {\it 
associative  comultiplication}. We  refer to (\ref{assiomap}) saying that 
$\Delta$ is {\it cocommutative}.

For any $X\in M$ we have
$$
 \Delta(X)= X\otimes 1+ 1\otimes X 
,\ \ \
\epsilon(X)=0
, \ \ \
\delta(X) = -X
.
$$ 
To give formulas for $\Delta$ we introduce the following notation.
Let $n\in\Naturali$ and $\Sigma_n$ be the group of permutations of $n$ 
elements. For any 
$s\in\Sigma_n$ and $X_1,.., X_n\in M $,   let 
$
   \alpha( X_{s(1)},...,X_{s(n)}) 
\in\{1, -1\} 
$ 
be   the sign such that
$
  \alpha( X_{s(1)},...,X_{s(n)})
  X_{s(1)}\cdot\cdot\cdot X_{s(n)} 
= X_1\cdot\cdot\cdot X_n$ in $S(M)
$. 
If $X\in M_0$
\begin{eqnarray}
\label{DeltaL}
  \Delta\left( X^n \right) 
= \sum_{j=0}^{n}
  {n\choose j}
  X^j\otimes X^{n-j}
, \ \forall n\geq 0
\end{eqnarray}
and, if $X_1,..., X_n\in M$ 
\begin{eqnarray*}
  \Delta\left( X_1\cdot\cdot\cdot X_n
        \right) 
&=& \sum_{j=0}^{n}
  \sum_{1\leq p_1<\cdot\cdot\cdot <p_j\leq n}
  \alpha(\vec X_{\vec p})
   X_{p_1}\cdot\cdot\cdot X_{p_j} \otimes 
   X_1\cdot\cdot\cdot\widehat{X_{p_1}}\cdot\cdot\cdot 
\widehat{X_{p_j}}\cdot\cdot\cdot X_n 
\end{eqnarray*}
where
$  \alpha(\vec X_{\vec p})
:= \alpha( X_{p_1},...,X_{p_j}
         , X_1,...,\widehat{X_{p_1}},...,\widehat{X_{p_j}},...,X_n
         )
$.

We denote $S(M)^*:= Hom(S(M), \campo)$. Because of $S(M)$ is a coalgebra, 
$S(M)^*$ is a commutative superalgebra and  it is called the algebra of 
formal power series over $M$. 

More generally, if $N$ is a $\campo$-module, $Hom( S(M), N)$ is called the 
space of {\it formal functions} on $M$ with values in $N$. 
Each $X\in N$ defines a "constant function" of  $Hom(S(M), N)$: it is the function such 
that  $1\mapsto X$ and $S^n(M)\mapsto \{0\}$ for $n\neq 0$. 
We have the following structure of  $S(M)^*$-module: 
$F\varphi  := (F\otimes\varphi)\circ\Delta$
for $\varphi\in S(M)^*$ and $F\in Hom(S(M), N)$.
Let $Y\in M$, we define 
\begin{eqnarray*}
\L(Y): Hom( S(M),N)\ni f 
\rightarrow (-1)^{p(f)p(Y)}f\circ Y^L\in Hom(S(M),N)
. 
\end{eqnarray*}
It is called the derivative in the direction $Y$. 

\begin{Oss}
\label{oss121}
By definition,  $\partial(Y)(X) = 0$ for any $X\in N$.
\end{Oss}

When $N=M$,   $Hom( S(M), M)$ is called the space of {\it formal vectors field} 
over $M$. 
The identity of $M$ extends to a morphism of $\campo$-modules 
$x_M:S(M)\rightarrow M$ by $S^n(M)\mapsto \{0\}$ for $n\neq 1$. 
It is called the {\it generic point of $M$}, and it will be denoted by $x$ 
when there is no risk of confusion.

\begin{Oss}
We have $\partial(Y)(x) = Y$ for any $Y\in M$.
\end{Oss}

Let $A$ be a   $\campo$-superalgebra. In $Hom(S(M), A)$ we have the  
following structure of $S(M)^*$-superalgebra: 
$     F\cdot G 
   := {\Mult}
\circ (F\otimes G)
\circ \Delta
$, 
for any $F,G\in Hom(S(M), A)$. 
\begin{Oss}
For any $Y\in M$, $\partial(Y)$ is a derivation of $Hom(S(M), A)$.
\end{Oss}
\noindent
We have seen that $A\subseteq Hom(S(M), A)$, moreover $A$ is a 
$\campo$-subsuperalgebra of 
$Hom(S(M), A)$.
If $A$ is associative  $Hom(S(M), A)$ is  associative, because $\Delta$ 
verifies  (\ref{lzero}).
If $A$ is unitary  $Hom(S(M), A)$ is  unitary, with unit given by 
$\epsilon : 
S(M)\ni W\mapsto 1\epsilon(W)\in A$.
If  $A$ is commutative $Hom(S(M), A)$ is commutative, because $\Delta$ is a 
cocommutative comultiplication. 

In the particular case $A =S(M)$, it is a tradition  to denote  by $*$ 
the product of $Hom(S(M), S(M))$. In this case $\delta\in Hom(S(M), S(M) )$ 
and identities (\ref{assiomaq}), (\ref{assiomaL}) give 
\begin{eqnarray}
\label{l10}
\delta*id=id*\delta =\epsilon  
.
\end{eqnarray}

\begin{lemma}
 If $\gG$ is a Lie $\campo$-superalgebra, $Hom(S(M),\gG)$ is a Lie 
$S(M)^*$-superalgebra. 
\end{lemma}

\subsection{Coderivations of a symmetric algebra}

Let $A$ be a $\campo$-module equipped of a comultiplication
$\Delta$. 

\begin{Def}
A coderivation of $A$ is a morphism of $\campo$-modules 
$\Phi :A\rightarrow A
$ 
such that 
$     \Delta\circ \Phi
    = (\Phi\otimes id + id\otimes \Phi)
\circ \Delta
$.
\end{Def} 

To describe the
coderivations of $S(M)$, we introduce the $\campo$-module
$  P(S(M)
 := \{ W\in S(M)| \ \Delta(W) = 1\otimes W + W\otimes 1
   \}
$. 
Its elements are called the {\it primitive elements of $S(M)$}. 
By definition of  $\Delta$, $M\subseteq P(S(M))$.

Let $\varphi : S(M)\rightarrow P(S(M))$ be a morphism of 
$\campo$-modules, we define  
$  \Phi 
:= id*\varphi   : S(M)\rightarrow S(M)
$.
\begin{tm}[Theorem 1, \cite{Rad}]
\label{lemma11}
\phantom{b}\ 
\newline
The map $\Phi$ is the unique a coderivation of $S(M)$ such that 
$
      \delta*\Phi
   =  \varphi
$.
\end{tm}
%\noindent \dimostrazione \ 
%As a first step  we show that $\Phi$ is a  coderivation. By 
%definition 
%\begin{eqnarray*}
%    &&(id\otimes \Phi)\circ\Delta 
%    = (id\otimes {\Mult})
%\circ (id\otimes id\otimes \varphi)
%\circ (id\otimes \Delta)
%\circ \Delta
%\\
%    &&(\Phi\otimes id)\circ\Delta 
%    = ({\Mult}\otimes id)
%\circ (id\otimes \varphi\otimes id)
%\circ (\Delta\otimes id)
%\circ \Delta
%\\
%    &&\Delta\circ \Phi
%    = \Big(
%              ({\Mult}\otimes id)
%        \circ (id\otimes \varphi\otimes id)
%        \circ (id\otimes \sigma)
%            + (id\otimes {\Mult}) 
%        \circ (id\otimes id\otimes \varphi)
%      \Big)
%\circ (\Delta\otimes id)
%\circ \Delta
%\end{eqnarray*}
%so identities (\ref{lzero}) and (\ref{assiomap}) give 
%that $\Phi$ is a coderivation. 
%
%The second part of the  theorem follows  from identities (\ref{l10}). In 
%fact  we get 
%$\varphi= \delta*id*\varphi = \delta*\Phi$.  
%\halmos 
%
%\begin{Oss}
%The previous theorem is valid if $S(M)$ is replaced by any 
%%$A$ which is a commutative superalgebra and a  
%cocommutative Hopf superalgebra. 
%\end{Oss}

\subsection{Generic point of a Lie superalgebra}

Let  $(\gG, [\cdot, \cdot])$ be a Lie $\campo$-superalgebra. For any 
$X\in \gG$, we denote by $\ad X$  the
application $[X,\cdot]:\gG\rightarrow \gG$.

Let $t$ and $u$ be two even commuting variables. For any 
$r, q\in\Naturali$ we introduce the notation 
\begin{eqnarray}
\label{not1}
   \left( t^r u^q        : [Y,Z] 
   \right)_X
:= [(\ad X)^r(Y), (\ad X)^q(Z)]
, \ \ \ 
\forall X\in\gG_0
, \ 
\forall Y,Z\in\gG
.
\end{eqnarray} 
By linearity it is extended to all polynomials in $\campo [t,u]$.

\begin{lemma}
\label{lemma12bis}
For all $q\in\campo[z]$, $X\in\gG_0$ and  $Y,Z\in\gG$, we have 
$$ 
  q(\ad  X)([Y,Z])
= \left( q(t+u): [Y,Z] \right)_X
%\in \gG
.
$$
\end{lemma}
\noindent\dimostrazione\
It is sufficient to consider $q(z) = z^k$ with $k\geq 1$. It means that 
it is sufficient to show that 
$
  (\ad X)^k([Y,Z])
= \sum_{p=0}^k 
 {k\choose p}
  [(\ad X)^p(Y), (\ad X)^{k-p}(Z)]
, \ k\geq 1
$. 
This identity means that $\ad X$ is an even derivation, which follows from the 
Jacobi identity. \ \halmos
\newline

We introduce $\gG_x := Hom(S(\gG),\gG)$.  Each $X\in\gG$ is 
identified to its image in  $\gG_x$. 
As seen above (section \ref{sec121}), the 
comultiplication $\Delta$ of $S(\gG)$ and the bracket for $\gG$ allow  to
define the bilinear application 
$      [F,G]
    := [\cdot,\cdot ] 
\circ (F\otimes G)
\circ\Delta
$, 
for any  $F,G\in\gG_x$. We have seen also  that $\gG_x$ is a 
Lie $S(\gG)^*$-superalgebra  and $\gG\subseteq\gG_x$ is a
Lie $\campo$-subsuperalgebra. Let $x\in\gG_x$ be the generic point of
$\gG$.

\begin{Oss}
\label{oss131}
For any  $n\in\Naturali\setminus\{0\}$, $(\ad
x)^n:\gG_x\rightarrow \gG_x$ is a 
$S(\gG)^*$-morphism. 
In particular, if  $Y\in\gG$, $(\ad x)^n(Y): S(\gG)\rightarrow\gG$ is the 
map  such that, for any $p\geq 0$ and $X_1,..., X_p\in\gG$ 
\begin{eqnarray*}
        X_1\cdot\cdot\cdot X_p
\mapsto 
\left\{
\begin{array}{ll}
        0
        ,& p\neq n

\\
        \sum_{s\in\Sigma_n} 
        \alpha(Y, X_s)
        \ad X_{s(1)}\circ\cdot\cdot\cdot\circ\ad X_{s(n)}
        (Y)
        ,& p=n
\end{array}
\right.
\end{eqnarray*}
where   
$
   \alpha(Y, X_s)
:= (-1)^{p(Y)p(X_1+\cdot\cdot\cdot +X_n)}
   \alpha( X_{s(1)},..., X_{s(n)} )
$. 

If $n=0$, 
$(\ad x)^0(Y) :=Y\in\gG_x$.
\end{Oss}
\noindent
Let $q = c_0+c_1t+c_2t^2+\cdot\cdot\cdot\in\campo[[t]]$  and $Y\in\gG$. 
As a consequence of  remark \ref{oss131}, we can define
$$  q(\ad x)(Y)
:= c_0 Y +c_1 (\ad x)(Y)+c_2(\ad x)^2(Y)+\cdot\cdot\cdot
.
$$
It is the morphism of  $\campo$-modules from $S(\gG)$ to $\gG$ such that,   
for any $n\in\Naturali$,   its restriction to $S^n(\gG)$ is $c_n(\ad x)^n(Y)$.

\begin{Oss}
\label{oss112}
In the Lie superalgebra $\gG_x$, we consider the formula 
(\ref{not1}) with  $X=x$. This gives the  formula of  a morphism of 
$\campo$-modules from $S(\gG)$ to $\gG$. Let 
$Y,Z\in\gG$, 
for any  $n\in\Naturali$ and $X_1,..., X_{n}\in\gG$, 
 this morphism   is  given by 
$$
\begin{array}{l}
(t^{r}u^{q} : [Y,Z])_x(X_1\cdot\cdot\cdot X_{n})=
\\
=\left\{
\begin{array}{ll}
0, & n\neq r+q
\\
        \sum_{s\in\Sigma_{n}} 
\alpha( X_{s})
[ \ad X_{s(1)}\circ\cdot\cdot\cdot\circ\ad X_{s(r)}(Y)
, \ad X_{s(r+1)}\circ\cdot\cdot\cdot\circ\ad X_{s(n)}(Z)
]
, &n=q+r
\end{array}
\right.
\end{array}
$$
where the coefficients $\alpha(X_s)$ are  given by
$$  \alpha(X_s)
:= (-1)^{ p(X_1+\cdot\cdot\cdot +X_n)p(Z) 
        + p(Y)p(X_{s(1)}+\cdot\cdot\cdot+ X_{s(r)})}
\alpha( X_{s(1)},..., X_{s(r+q)} )
.
$$
 
\end{Oss}
\noindent
As above this allows to define $(p(t,u):[Y,Z])_x$ for any 
formal power 
series 
$p\in\campo[[t,u]]$. 
The following theorem plays a basic role in this text.

\begin{tm}
\label{cor11}
Let $Y,Z\in\gG$ and $q(z)\in\campo[[z]]$. In $\gG_x$ we have    
$$
    \L(Y)\big(  q(\ad x)(Z) \big)
  = (-1)^{p(q)p(Y)}
    \left(
        \frac{ q(t+u) - q(u) }{t}
        : [Y,Z]
    \right)
.
$$
\end{tm}
\noindent\dimostrazione \ 
We only  need to consider the case $q(z) = z^k$,  with  $k\geq 0$. For $k=0$ the 
statement follows from remark \ref{oss121}. We recall that $\L(Y)$ is a 
derivation. By 
induction over $k$ and by the Jacobi identity
in $\gG_x$, we get
\begin{eqnarray*}
&&  \L(Y)( (\ad x)^{k+1}(Z) ) 
=   \L(Y)( [x, (\ad x)^{k}(Z)] ) =
\\
&=& [Y, (\ad x)^k (Z)]
  + [x, \L(Y)((\ad x)^k)(Z)]
\\
&=&  ( u^k : [ Y,  Z])
  + \ad x
    \left( \left(\frac{(u+t)^k - u^k}{t} : [Y, Z] \right)_x
    \right)
\\
&=& ( u^k : [ Y,  Z])_x
  + \left( (t+u)\frac{(u+t)^k - u^k}{t} : [Y, Z]
    \right)_x
\\
&=& \left( \frac{(u+t)^{k+1} - u^{k+1}}{t} : [Y, Z]
    \right)_x
.\  \halmos
\end{eqnarray*}

\section{Functional equations  associated to coderivations}
\label{sec21}

Let $\campo$ be a  commutative superring and 
$\varphi(z) = \sum_j z^jc_j \in \campo[[z]]$.
For  any  Lie $\campo$-superalgebra $\gG$ and  $a\in\gG$, we define the 
formal vector field
\begin{equation}
\label{secondeL}
 \varphi^a
:=\varphi(\ad x)(a)
\in \gG_x
.  
\end{equation}
We recall from remark  (\ref{oss131}) that  
$$
  \varphi^a(X_1\cdot\cdot\cdot X_n)
= (-1)^{p(a)p(X_1+\cdot\cdot\cdot +X_n)}
  \sum_{\sigma\in\Sigma_n}
  c_n
  \alpha(X_{\sigma(1)},..., X_{\sigma(n)})
      \ad X_{\sigma(1)}
\circ \cdot\cdot\cdot
\circ \ad X_{\sigma(n)}
      (a)
$$
for any $n\in\Naturali$ and $X_1,..., X_n\in\gG$. 
In particular, if $X\in\gG_0$ we get 
\begin{eqnarray}
\label{phi1}
&&   \varphi^a\left( 1\right)
 =    c_0 a
\\\nonumber
&&   \varphi^a\left( X^n \right)
 = n!
   c_n
   (\ad X)^n
   (a)
, \ \ \forall n\geq 1
.
\end{eqnarray}

\begin{Oss} {\bf (Functorial property)}
\label{oss11}
\newline
Let $\hG$ be a Lie $\campo$-superalgebra and $f:\gG\rightarrow \hG$ be a 
morphism of Lie $\campo$-superalgebras. 
The formula (\ref{secondeL}) shows that
$
  f\circ \varphi^a
  (X_1\cdot\cdot\cdot X_n) 
= \varphi^{f(a)}
  (f(X_1)\cdot\cdot\cdot f(X_n)) 
$, 
for any $n\in\Naturali$ and $X_1,..., X_n\in\gG$.
\end{Oss}
\noindent

Let 
$
  \rho(t,u) 
= \sum_{j\geq 0}
  \sum_{i=0}^j
  t^iu^{j-i}
  d_{i,j-i} 
\in\campo[[t,u]]
$. 
To  $a,\ b\in\gG$ we  associate  the formal vector field  
$(\rho(t, u): [a,b])_x\in\gG_x$.  
We recall that (remark \ref{oss112}) for any $X\in\gG_0$  and 
$n\in\Naturali$ 
we get 
$
  (\rho(t,u): [a,b])_x (X^n)
= n!
  \sum_{i=0}^n
  \left(  t^iu^{n-i}
          d_{i, n-i}
        : [a, b]
  \right)_X
$.

\begin{Oss}
\label{oss211}
By lemma \ref{lemma12bis} we have 
$ \varphi^{[a,b]} 
= ( \varphi(t+u):[a,b])_x
$. 
\end{Oss}

\noindent 
By theorem \ref{lemma11}, to the formal vector field $\varphi^a$ we 
associate  the
coderivation 
\begin{equation}
\label{terzaL}
       \Phi^a
    := id*\varphi^a
\equiv {\Mult}\circ (1\otimes \varphi^a)\circ\Delta
. 
\end{equation}
For any $n\in\Naturali$ and $X_1,..., X_n\in\gG$ we have 
$$
  \Phi^a
  \left(  X_1\cdot\cdot\cdot X_n
  \right) 
= \sum_{j=0}^{n}
  \sum_{1\leq p_1<\cdot\cdot\cdot <p_j\leq n}
  a(\vec X, \vec p) 
  X_{p_1}\cdot\cdot\cdot X_{p_j}
\cdot\varphi^a
   \left(
   X_1\cdot\cdot\cdot\widehat{X_{p_1}}\cdot\cdot\cdot\widehat{X_{p_j}}
\cdot\cdot\cdot X_n 
   \right)
,
$$
where 
$
   a(\vec X, \vec p) 
:= (-1)^{p(\varphi^a)p(X_{p_1}+\cdot\cdot\cdot+X_{p_j})}
   \alpha(    X_{p_1},..., X_{p_j}, X_1, 
         ..., \widehat{X_{p_1}},..., \widehat{X_{p_j}},...,  X_n
         )
$. 
In particular,  if $X\in\gG_0$ we get 
$
     \Phi^a\left(X^n \right)
   =  \sum_{j=0}^n
      {n\choose j}
      X^j\cdot\varphi^a(X^{n-j})
$. 

Let $\psi\in\campo[[t]]$, $b\in\gG$, and let  
$\Psi^b :S(\gG)\rightarrow S(\gG)$ be the associated coderivation. 

\begin{Oss}
\label{rem213}
By definition, for any $Y\in\gG$ we have  $id* Y = Y^L$.
\end{Oss}

\begin{lemma}
\label{lemma21}
For any  $Y\in\gG$ we have 
\newline
i)
$ 
      \Phi^a\circ Y^L
    = id
    * \left(
         \varphi^a* Y
        -\left(
          \frac{\varphi(t+u) - \varphi(t)}{u}
        : [a,Y]
          \right)_{x}
\right)
$
\newline
\newline
ii) 
$
      \Phi^a\circ \Psi^b
    = id
    * \left(
          \varphi^a* \psi^b
        - (-1)^{p(a)p(\psi)}
          \left(
          \frac{\varphi(t+u) - \varphi(t)}{u}
          \psi(u)
        : [a,b]
          \right)_{x}
\right)
$.
\end{lemma}
\noindent\dimostrazione \ 
$i)$
From the fact that $\Delta$ is, in particular,  a morphism of algebras, 
and from the remark \ref{rem213} we have 
$
  \Phi^a\circ Y^L
  \equiv \big( 
        id*\varphi^a 
    \big)
\circ Y^L
    = id
    * \left\{
          \varphi^a\circ Y^L
    + (-1)^{p(\varphi^a)p(Y)}
       Y
    * \varphi^a
      \right\}
$. 
As $*$ is commutative,  this shows that 
$
    \Phi^a\circ Y^L
  = id
  * \left\{
          \varphi^a\circ Y^L
      \right\}
  + id
  * \varphi^a*Y
$. 
By definition 
$
  \varphi^a\circ Y^L
= (-1)^{p(Y)(p(\varphi)+p(a))}
  \L(Y)(\varphi^a)
$, 
so the theorem \ref{cor11} gives the desired formula.

$ii)$ 
Let us consider the Lie superalgebra $\gG_x$ and its generic point 
$y\in Hom( S(\gG_x), \gG_x)$. If $X_1,..., X_n\in\gG\subset\gG_x$ we 
have 
$
  \left(   \Phi^a_\gG\circ \Psi^b_\gG
  \right)
  (  X_1\cdot\cdot\cdot X_n)
= 
  \left(   \Phi^a_{\gG_x}\circ \Psi^b_{\gG_x}
  \right)
  (  X_1\cdot\cdot\cdot X_n)
$, so it is sufficient to prove the statement for $\gG_x$. 
By definition and by remark \ref{rem213} 
we have 
\begin{eqnarray*}
    \Phi^a_{\gG_x}\circ \Psi^b_{\gG_x}
   =  \Phi^a_{\gG_x}
\circ \left( id*\psi(\ad y)(b) 
      \right)
    = \Phi^a_{\gG_x}
\circ \psi(\ad y)(b)^L
,
\end{eqnarray*}
so from case $i$, we get 
\begin{eqnarray*}
    \Phi^a_{\gG_x}\circ \Psi^b_{\gG_x}
&=& id
  * \left\{
          \varphi(\ad y)(a)
        * \psi(\ad y)(b)
        - \left(
                  \frac{\varphi(t+u) - \varphi(t)}{u}
                : [a, \psi(\ad y)(b)]
          \right)_{y}
\right\}
. 
\end{eqnarray*}
By definition 
$  \left(
          \frac{\varphi(t+u) - \varphi(t)}{u}
        : [a, \psi(\ad y)(b)]
  \right)_{y}
= (-1)^{p(a)p(\psi)}
  \left(
          \frac{\varphi(t+u) - \varphi(t)}{u}
          \psi(t)
        : [ a,b ]
  \right)_{y}
$ 
so the proof is finished.\ \halmos

\begin{tm}
\label{cor21}
\begin{eqnarray*}
  [\Phi^a, \Psi^b]
= id
* (-1)^{p(\psi)p(a)}
  \left(
       - \frac{\varphi(t+u) - \varphi(t)}{u}
         \psi(u)
       - \varphi(t)\frac{\psi(t+u) - \psi(u)}{t}
        : [a,b]
  \right)_x
.
\end{eqnarray*}
\end{tm}
\noindent\dimostrazione\ 
Let 
$ 
  \omega(t,u) 
:= -\frac{\varphi(t+u) - \varphi(t)}{u}
   \psi(u)
 - \varphi(t)\frac{\psi(t+u) - \psi(u)}{t}
$, 
 we denote by $\Omega^{[a,b]}$ the coderivation
corresponding to 
$(-1)^{p(\psi)p(a)}
 \left(
        \omega(t,u) : [a,b]
 \right)_x
$. 
By theorem 
\ref{lemma11} we want to 
show that 
$[\Phi^a, \Psi^b]
=\Omega^{[a,b]}
$. 
From lemma \ref{lemma21} we get 
\begin{eqnarray*}
    &&  [\Phi^a, \Psi^b]=
\\
  &=& id{*}
\left(
        - (-1)^{p(a)p(\psi)}
          \left(
          \frac{\varphi(t+u) - \varphi(t)}{u}
          \psi(u)
        : [a,b]
          \right)_{x}
        + \varphi^a*\psi^b
\right)+
\\
&& id{*}
\left(
(-1)^{p(b)p(\varphi)+p(\Phi^a)p(\Psi^b)}
          \left(
          \frac{\psi(t+u) - \psi(t)}{u}
          \varphi(u)
        : [b,a]
          \right)_{x}
        -(-1)^{p(\Phi^a)p(\Psi^b)}
         \psi^b*\varphi^a
\right)
\\
&=& - (-1)^{p(a)p(\psi)}
   id* \left(
          (-1)^{        p(\varphi)p(\psi)}
          \frac{\psi(t+u) - \psi(u)}{t}
          \varphi(t)
        + \frac{\varphi(t+u) - \varphi(t)}{u}
          \psi(u)   
        : [a,b]
     \right)_{x}
\\
&=& id*  (-1)^{p(\psi)p(a)}
    \left( \omega(t,u): [a,b]
    \right)_{x}
. \ \halmos
\end{eqnarray*}
\newline
To prove the next theorem we need some preliminaries, which we 
state  in a form that will be useful later.

\begin{Def}
Let $N\geq 1$ be an integer. A Lie superalgebra $\gG$ is said to be 
$N$-nilpotent if we have 
$
\ad X_1\circ \cdot\cdot\cdot\circ \ad X_{N}=0
$ for any $X_1,..., X_N\in\gG$.
\end{Def}
\begin{Oss}
For $N=1$ we have a commutative Lie superalgebra, for $N= 2$  we have  
a Lie superalgebra of Heisenberg type. 
\end{Oss}

\begin{lemma}
\label{lem212}
For any  $N\geq 2$, there exists a $N$-nilpotent Lie $\campo$-superalgebra 
$\gG_N$, equipped of an infinite family of even elements 
$\{\alpha,\beta, X_1, X_2,...\}$
such that 
$$
\bigcup_{0\leq r+s\leq N-2}
\{ [ \ad X_{i(1)}\circ\cdot\cdot\cdot\circ \ad X_{i(r)}(\alpha)
   , \ad X_{i(r+1)}\circ\cdot\cdot\cdot\circ \ad X_{i(r+s)}(\beta)
   ]
 ; i(1),..., i(r+s)\in\Naturali
\}
$$ 
is a contained in a basis. 
\end{lemma}
\noindent\dimostrazione\ 
We start by considering the free Lie algebras $\hG$ over $\Z$, with an 
infinite family of  generators $\alpha$, 
$\beta$, $X_1,X_2,..$. By properties of free Lie algebras 
({\bf\cite{Bou2}}, prop. 10, page 
26) we know that $\hG$ is free, and that 
$$
\bigcup_{r,s\geq 0}
\{ [ \ad X_{i(1)}\circ\cdot\cdot\cdot\circ \ad X_{i(r)}(\alpha)
   , \ad X_{i(r+1)}\circ\cdot\cdot\cdot\circ \ad X_{i(r+s)}(\beta)
   ]
 ; i(1),..., i(r+s)\in\Naturali
\}
$$ 
is contained in a basis of $\hG$.

Let $I_N$ be the ideal of $\hG$ 
generated by 
$\{ \ad x_1\circ\cdot\cdot\cdot\circ\ad x_N (Y)| x_1,..., x_N, Y\in\hG\}$. 
The quotient $\hG_N:=\hG/I_N$ is a $N$-nilpotent Lie superalgebra 
over $\Z$ and the family 
$
f_N:= \bigcup_{0\leq r+s\leq N-2}
\{ [ \ad X_{i(1)}\circ\cdot\cdot\cdot\circ \ad X_{i(r)}(\alpha)
   , \ad X_{i(r+1)}\circ\cdot\cdot\cdot\circ \ad X_{i(r+s)}(\beta)
   ]
 ; i(1),..., i(r+s)\in\Naturali
\}
$ is contained in a basis of $\hG_N$.  

We define $\gG_N:= \hG_N\otimes \campo$. It is a 
$N$-nilpotent  Lie superalgebra over 
$\campo$ and 
% and $\hG_N$ is a Lie subsuperalgebra. 
%As the tensor product of modules is distributive, $\gG_N$  inherits from 
%$\hG_N$ the property that 
$f_N$ is contained in a basis of $\gG_N$.  
\halmos

\begin{lemma}
\label{lem213}
Let 
$  \omega(t,u)
 = \sum_{i,j=0}^\infty
   c_{ij}
   t^iu^j
\in\campo[[t,u]]
$ and $N\geq 2$.  
If  for any $N$-nilpotent Lie $\campo$-superalgebra $\gG$ we have 
$
  (\omega(t,u): [a,b])_x
= 0
, \ \forall a,b\in\gG
$,  
then $c_{ij}=0$ for any $0\leq i+j\leq 
N-2$.
\end{lemma}
\noindent\dimostrazione\  
 We consider the case $\gG=\gG_N$, where $\gG_N$ is the $N$-nilpotent Lie 
superalgebra of lemma \ref{lem212}.
%If $\gG=\gG_2$ we  have $0 = (\omega(t,u): [a,b])_x = 
%c_{00}[a,b]$. 
Choosing $a=\alpha$ and $b=\beta$ 
we get 
$$  (\omega(t,u), [\alpha,\beta])_x
 = \sum_{i+j=0}^{N-2}
   c_{i,j}
   (t^iu^{j}:   [\alpha,\beta])_x
.
$$ 
Let $0\leq p\leq N-2$, the 
remark \ref{oss112} gives 
\begin{eqnarray*}
&&(\omega(t,u), [\alpha,\beta])_x(X_1\cdot\cdot\cdot X_p)=
\\
&=&  \sum_{i=0}^{p}
   c_{i,p-i}
   \sum_{s\in\Sigma_p}
   [ \ad X_{s(1)}\circ\cdot\cdot\cdot \circ\ad X_{s(i)}(\alpha),
   , \ad X_{s(i+1)}\circ\cdot\cdot\cdot \circ\ad X_{s(p)}(\beta)
   ]
.
\end{eqnarray*} 
As $(\omega(t,u), [\alpha,\beta])_x(X_1\cdot\cdot\cdot X_p)$ is zero, 
lemma \ref{lem212} gives that    $c_{i,p-i}=0$ for any $i=0,..., p$. As  
$0\leq p\leq N-2$, the proof is finished. \halmos
\newline

Let $\lambda\in\campo[[z]]$. For any $a\in\gG$ we consider the 
the formal vector field $\lambda^a$ (see formula \ref{secondeL}) 
and the corresponding coderivation   $\Lambda^{a} : S(\gG)\rightarrow 
S(\gG)$ (see formula \ref{terzaL}).

\begin{tm}
\label{prop21}
Let $\varphi, \psi, \lambda\in\campo_0[[t]]$. 
For any Lie $\campo$-superalgebra $\gG$ we have
\begin{eqnarray}
\label{doppioAst}
  [\Phi^a, \Psi^b] 
= \Lambda^{[a,b]}
, \ \ \forall a, b\in\gG
\end{eqnarray}
if and only if 
$
    \varphi
  , \psi
  , \lambda
$ 
verify
\begin{eqnarray*}
\left(
- \frac{\varphi(t+u) - \varphi(t)}
       {u}
  \psi(u)
- \varphi(t)
  \frac{\psi(t+u) - \psi(u)}
       {t}
\right)
= \lambda(t+u)
\end{eqnarray*}
in $\campo_0[[t,u]]$.
\end{tm}
\noindent\dimostrazione\ 
Let 
$
   \omega(t,u)
:= - \frac{\varphi(t+u) - \varphi(t)}
       {u}
  \psi(u)
- \varphi(t)
  \frac{\psi(t+u) - \psi(u)}
       {t}
- \lambda(t+u)
$. 
Using theorem \ref{cor21} and  remark 
\ref{oss211}, we see that  (\ref{doppioAst}) is equivalent to 
$
  id
* \left(
    \omega(t,u)
  : [a,b]
  \right)_x
= 0
, \forall a,b\in\gG
$. 
By theorem \ref{lemma11},  this identity is equivalent to 
$
 \left(
    \omega(t,u)
  : [a,b]
  \right)_x
= 0
, \forall a,b\in\gG
$. 

We get immediately that the functional equation 
is sufficient. To show the converse, it is sufficient to apply the lemma 
\ref{lem213} to any $N$-nilpotent Lie superalgebra $\gG_N$, with $N\geq 
2$. We get that in $\omega(t,u)$ the coefficients of degree $N-2$ are 
zero, for any $N\geq 2$. In particular $\omega(t,u)=0$. \halmos

\begin{tm}
\label{cor31}
Let $\varphi\in\campo[[t]]$. For any Lie $\campo$-superalgebra $\gG$, 
we have
\begin{eqnarray}
\label{eq26}
  [\Phi^a, \Phi^b] 
= \Phi^{[a,b]}
, \ \ \ \forall a, b\in\gG
\end{eqnarray}
if and only if $\varphi$  has  
even coefficients ($\varphi\in \campo_0[[t]]$) and  verifies 
\begin{eqnarray}
\label{precL}
%\label{eq15}
\left(
- \frac{\varphi(t+u) - \varphi(t)}
       {u}
  \varphi(u)
- \varphi(t)
  \frac{\varphi(t+u) - \varphi(u)}
       {t}
\right)
= \varphi(t+u)
.
\end{eqnarray}
\end{tm}
\noindent\dimostrazione\ 
As $p(\Phi^a) \equiv p(\varphi)+p(a)$, the identity  (\ref{eq26}) needs  
$p(\varphi)=0$. The theorem \ref{cor31} follows from theorem \ref{prop21}. 
\halmos

\section{Universal representations}
\label{sec222} 
Let $\campo$ be a commutative superring,   
$ \varphi(t)
= \sum_j t^jc_j
  \in \campo_0[[t]]
$, and 
$\gG$ be a   Lie $\campo$-superalgebra. We consider the map 
$\Phi:\gG\ni a\mapsto \Phi^a\in Hom(S(\gG), S(\gG))$ defined in 
(\ref{terzaL}).
From theorem \ref{cor31} we know  that  $\Phi$ is a
representation for any Lie $\campo$-superalgebras $\gG$,  if and
only 
if $\varphi\in\campo_0[[t]]$ verifies the functional 
equation (\ref{precL}). 
Before looking for solutions of this functional equation,  we 
introduce the notion of {\it universal representation}. 

Let $M(\gG)$ be the symmetric algebra or the enveloping algebra 
of $\gG$ (defined in paragraph \ref{sec23}).

\begin{Def}
\label{def32}
Assume that 
\newline
i) for any $\campo$-Lie superalgebra $\gG$, we have a
representation $\Phi_\gG:\gG\rightarrow End(M(\gG))$. 
\newline 
ii)  For any $a\in\gG$ and any morphism of Lie 
$\campo$-superalgebras $f:\gG\rightarrow  \hG$ the 
following diagram, with $\tilde f: M(\gG)\rightarrow M(\hG)$ the 
induced morphism of algebras, is commutative:
$$
\begin{array}{ccc}
M(\gG)&\stackrel{\Phi^a_\gG}{\longrightarrow} &M(\gG)
\\
\downarrow_{\tilde f} &&\downarrow_{\tilde f}
\\
M(\hG)&\stackrel{ \Phi^{f(a)}_\hG }{\longrightarrow} &M(\hG)
\end{array}
.
$$
Then we say that  $\Phi$ is  an universal representation in the category 
of Lie $\campo$-superalgebras.

Let $N\in\{2,3,4,...\}$. In an analogous way we define the universal 
representations in the category 
of 
$N$-nilpotent Lie $\campo$-superalgebras.
\end{Def}

For any commutative superring $\campo$ we introduce   
\begin{eqnarray}
\varphi_{0}(t) := -t
\in\campo_0[[t]]
.
\end{eqnarray}
For any $c\in\campo^\times_0$, if 
$\campo\supseteq\Razionali$, we  introduce also 
\begin{eqnarray}
  \varphi_{c}(t)
= \frac{t}{e^{\frac{t}{c} }-1}
\in\campo_0[[t]]
.
\end{eqnarray}
All these 
series verify 
 $\varphi_c(0) = c$.

\begin{lemma}
\label{lemma23bis}
Let $\campo_0$ be a commutative field. The      solutions 
of  equation (\ref{precL}) which lie in $\campo_0[[t]]$ and such that 
$\varphi(0)=0$, are $\varphi = 0$ and $\varphi= \varphi_0$. 
\end{lemma}
\noindent\dimostrazione \ 
If $\varphi(0)=0 $, the limit $\lim_{u\rightarrow 0}$ applied to the 
equation (\ref{precL}) gives 
$
  \varphi(t)
  \left( 1 + \frac{\varphi(t)}{t} \right)
= 0
$, 
so $\varphi(t)$ is zero or $\varphi(t) = -t$ because $\campo_0[[t]]$ 
is 
a domain. \halmos

\begin{tm}
\label{lemma23}
i) Let  $\varphi\in\campo_0[[t]]$ be a solution of
equation (\ref{precL}). If the constant term   $\varphi(0)=:c$ is 
invertible, then 
$  f(t)
:= \frac{\varphi(t)+t}{\varphi(t)}
$  satisfies 
\begin{eqnarray}
\label{expoL}
\left\{
\begin{array}{l}
  f(t)\cdot f(u)
= f(t+u)
\\
f(0)=1
\\
f'(0)=\frac{1}{c}
\end{array}
\right.
.
\end{eqnarray}
ii) The system (\ref{expoL}) has  solutions if and only if  $\campo_0$ 
 contains $\Razionali$. In this case the unique solution is 
$e^{\frac{t}{c}}\in\campo_0[[t]]$. 
\newline
iii) Let $\campo\supseteq\Razionali$ and $c \in\campo^\times_0$.
 The unique solution of 
(\ref{precL}) 
in $\campo_0[[t]]$ verifying
  $\varphi(0)=c$ is $\varphi_c(t)$.
\end{tm}
\noindent\dimostrazione\ 
$i)$  We recall that $c$ is invertible if and only if    the  
series 
$\varphi$ is 
invertible, so we  
write the equation (\ref{precL}) as 
$$
     \frac{\varphi(t)+t}{\varphi(t)}
\cdot\frac{\varphi(u)+u}{\varphi(u)}
   = \frac{\varphi(t+u)+t+u}{\varphi(t+u)}
.
$$
We have 
$ \frac{\varphi(t)+t}{\varphi(t)}
= 1+\frac{1}{c}t+\cdot\cdot\cdot
$.
\newline
$ii)$ Let $f=1+\frac{1}{c}t+\sum_{k=2}^\infty f_kt^k$. The system 
(\ref{expoL}) gives $f'(t)= \frac{1}{c}f(t)$, so  
$2f_2=\frac{1}{c^2}$ 
and 
$kf_k=\frac{1}{c}f_{k-1}$ for any $k\geq 3$. 
By induction we get that $k$ is invertible and 
$f_k= \frac{1}{k!c^k}$ for 
any $k\geq 2$. 
\newline
$iii)$ When 
$ 
  f(t)
= e^{\frac{1}{c}t}
$,  we get  $\varphi(t) = \varphi_c(t)$. \ \halmos

\begin{Oss}
The Bernoulli numbers $\{b_k\in\Razionali, k\in\Naturali\}$ are defined by 
the generating series 
$$
  \varphi_{1}(z)
\equiv \frac{1}{e^z-1}
= 
  \sum_{k\geq 0}^{\infty}
  \frac{b_k}{k!}z^k
.
$$
For example $b_0=1$, $b_1= -\frac{1}{2}$, $b_2= \frac{1}{6}$. 
Let $c\in\campo^\times_0$, the fact that 
$\varphi_c(t)\in\campo_0[[t]]$ verifies the identity (\ref{precL})
can be written in the following way:
$$  
  \forall k\geq 0
, \ \ 
  0 
= b_k
+ \sum_{p=0}^{k-l}
  {k-l\choose p}
  \frac{b_p b_{k+1-p} }{ (l+1) }
+ \sum_{p=0}^{l}
  {l\choose p}
  \frac{b_p b_{k+1-p} }
  {(k+1-l)}
, \ \ \ l=0,..., k
.
$$
\end{Oss}

We have shown the following theorems

\begin{tm}
\label{tm213}
The map $\Phi_0 : \gG\rightarrow Hom( S(\gG), S(\gG))$  
associated to $\varphi_0$ is a
representation by coderivations. 
\end{tm}

\begin{Oss} 
%\label{219}
Let $a\in\gG$. The map $\Phi_0^a$ is in the same time a derivation 
and a coderivation of  $S(\gG)$ : for any $X_1,..., X_n\in\gG$ we have 
$$
  \Phi_0^a(X_1\cdot\cdot\cdot X_n) 
= \sum_{j=1}^n
  (-1)^{p(a)p(X_1+\cdot\cdot\cdot + X_{j-1})}
  X_1\cdot\cdot\cdot \Phi_0^a(X_j)\cdot\cdot\cdot X_n
.
$$
It is the only derivation of $S(\gG)$ such that $\Phi^a_0(X)= [a, X]$ 
for $X\in\gG$, so  $\Phi_0$ is the adjoint representation of $\gG$ in 
 $S(\gG)$.
\end{Oss}
 
\begin{tm}
\label{tm222}
Let $\campo\supseteq\Razionali$. For any $c\in\campo^\times_0 
$, 
the series
$
\varphi_c(z)\in\campo_0[[z]]
$ 
gives a representation by coderivations 
$\Phi_c : \gG\rightarrow Hom( S(\gG), S(\gG) )$. 
\end{tm}

Let $\gG$ and $\hG$ be two Lie $\campo$-superalgebras,  $f:\gG\rightarrow 
\hG$ be a morphism of Lie $\campo$-superalgebras. It extends to a morphism 
of $\campo$-superalgebras $\tilde f : S(\gG)\rightarrow S(\hG).$

\begin{Oss} {\bf (Functorial property)}
\label{rem219}
\newline
By remark \ref{oss11}, for any $a\in\gG$ and $c\in\campo^\times_0$, the 
following diagram 
commutes
$$
\begin{array}{ccc}
S(\gG)&\stackrel{\Phi_c^a}{\longrightarrow} &S(\gG)
\\
\downarrow_{\tilde f} &&\downarrow_{\tilde f}
\\
S(\hG)&\stackrel{ \Phi_c^{f(a)} }{\longrightarrow} &S(\hG)
\end{array}
.
$$
In particular, $\Phi_0$ and 
$\Phi_c$ with $c\in\campo^{\times}_0$ are 
universal representations by coderivations.
\end{Oss}
\subsection{The case of nilpotent Lie superalgebras}
\label{subsec221}
We  give an analogue of theorem \ref{tm222} for $\campo$ 
not necessarily containing $\Razionali$. 
%We start recalling the following 
%definition.
%
%

Let $N\geq 2$,  $\gG$ a  $N$-nilpotent Lie superalgebra over $\campo$, 
$a$ and $b\in\gG$.

\begin{Oss}
i) The notation $\varphi(\ad x)(a)\in\gG_x $ is well-defined for 
$\varphi\in\campo[t]/t^N$ a truncated polynomial with coefficients 
in $\campo$.
\newline
ii) The notation $ (\rho(t,u) : [a,b])_x\in\gG_x$ is well-defined 
if $\rho(t,u)\in\campo[t,u]/I_N$, where $I_N$ is the ideal generated by 
$\{t^iu^j, i+j\geq N-1\}$. 
\end{Oss}
\noindent
To a truncated polynomial $\varphi(z)\in\campo_0[t]/t^{N}$, 
 we associate by formulas (\ref{secondeL}) and (\ref{terzaL}),
a family of coderivations still denoted by $\Phi^a:S(\gG)\rightarrow 
S(\gG)$,
$a\in\gG$.

\begin{tm}
For any $\gG$ a $N$-nilpotent  Lie 
superalgebra over $\campo$ the map $\Phi:\gG\ni a\mapsto \Phi^a\in 
Hom(S(\gG), S(\gG))$ is a 
representation by coderivations, if and only if 
$\varphi$   
verifies 
\begin{eqnarray}
\label{la217}
  \varphi(u)
  \frac{ \varphi(t+u)- \varphi(t)}
  {u}
+ \varphi(t)
  \frac{ \varphi(t+u)- \varphi(u)}
  {t}
= -\varphi(t+u)
\end{eqnarray}
in $\campo_0[t,u]/I_N$.
\end{tm}
\noindent\dimostrazione\ 
The direct part of the following theorem is a particular case 
of theorem 
\ref{prop21}. Let us prove the converse. 
Let 
$
   \omega(t,u)
:= -\varphi(u)
   \frac{ \varphi(t+u)- \varphi(t)}
   {u}
 - \varphi(t)
   \frac{ \varphi(t+u)- \varphi(u)}
   {t}
 - \varphi(t+u)
$. 
Proceeding as in the proof of theorem \ref{prop21} we get that $\Phi$ is a 
representation if and only if 
$
  \left(  \omega(t,u) : [a,b]
  \right)_x
= 0
, \ \forall a,b\in\gG
$. 
Moreover, for a $N$-nilpotent Lie superalgebra, this reduces to  
$$
  \left(  \omega(t,u) \mod I_N: [a,b]
  \right)_x
= 0
, \ \forall a,b\in\gG
.
$$ 
Using lemma \ref{lem213} we see that 
$ \omega(t,u) \mod I_N = 0 $. \halmos

\begin{ex}
Let $\campo_0$ be a field.
\newline 
i) Let $N=2$ and $\frac{1}{2}\in\campo_0$. We look for  $\varphi(t) = c_0+ 
c_1t \mod t^2$ solution of 
$
  2c_0c_1
= -c_0
$. 
We get $\varphi(t)= c_1t$ or $\varphi(t)=c_0-\frac{1}{2}t$.
\newline 
ii) Let  $N=3$ and $\frac{1}{2}, \frac{1}{3}\in\campo_0$. We look for 
$\varphi(t) = c_0+ 
c_1t+c_2t^2 \mod t^3$ 
solution of  
$
  2c_0c_1
+ (3c_0c_2+c_1^2)
  (u+t)
= -c_0
- c_1(t+u)
$. 
We get 
$\varphi(t)=c_2t^2$, or $\varphi(t)=-t+c_2t^2$, or $\varphi(t)=c_0 
-\frac{1}{2}t +\frac{1}{12c_0}t^2$ with $c_0\neq 0$.
\end{ex}

\begin{lemma}
Let $N\geq 2$.   The equation (\ref{la217}) has solutions in 
$\campo[t]/t^N$ with $\varphi(0)\in\campo^\times_0$,  if and only if 
 $\frac{1}{2},..., \frac{1}{N}\in\campo$. In this case the unique  
solution such that $\varphi(0)=:c\in\campo^\times_0$ is  
$$
        \varphi_{c}(t)
        \ mod \ t^{N}
. 
$$ 
\end{lemma}
\noindent\dimostrazione\
We look for $\varphi\in\campo_0[t]/t^{N}$ such that  
$1+\frac{t}{\varphi(t)}\in\campo_0[t]/t^{N+1}$  solves the system
(\ref{expoL}) in $\campo[t,u]/I_{N+2}$. 

The system (\ref{expoL}) has solutions in 
$\campo_0[t,u]/I_{N+2}$ exactly when   $2,...,N$ are invertible in 
$\campo$. In this case the unique solution is  $e^{\frac{1}{c}t} 
\mod t^{N+1}$ with 
$c\in\campo^\times_0$. It means that $\varphi(t) = \varphi_c(t)$. 
\halmos
\newline
\newline
We have shown that

\begin{tm}
\label{tm216}
Let $\campo\ni\frac{1}{2},...,\frac{1}{N}$. 
For any $c\in\campo^\times_0 $, the truncated polynomial  
$
\varphi_c(t)
\in\campo_0[t]/t^n
$ 
gives a representation by coderivations 
$\Phi_c : \gG\rightarrow Hom( S(\gG), S(\gG) )$. 
\end{tm}

\begin{Oss}
Let $\campo$ be a field, $p$ its characteristic. 
 When $2\leq N<p$, the previous theorem applies.
\end{Oss}

\subsection{Some properties of the representations $\Phi_c$}

This section applies to the  case $\campo\supseteq\Razionali$ and $\gG$ 
any Lie  $\campo$-superalgebra, and to the case 
$\campo\ni \frac{1}{2},..., \frac{1}{N}$ with $N\geq 2$ and $\gG$ a 
$N$-nilpotent Lie $\campo$-superalgebra.

\begin{Oss} 
If $c\in\campo^\times_0$, the representation 
$\gG\ni a\mapsto \Phi^a_c\in Hom(S(\gG), S(\gG))$ is faithful because 
$\Phi^a_c(1) = c \cdot a$.
\end{Oss}

\begin{Oss}
\label{oss217}
The  theorem \ref{prop21} gives that   $[\Phi_0^a,\Phi_c^b]
=\Phi_c^{[a,b]}$ for any  $a,b\in\gG$,  $c\in\campo^\times_0\cup\{0\}$.
\end{Oss}

\begin{tm}
\label{cor213}
Each  representation $\Phi_c$, $c\in\campo^\times_0$, is equivalent to
$\Phi_1$.
\end{tm}
\noindent\dimostrazione\ 
Let  $c\in\campo^\times_0$. We consider the map  $f_c: 
S(\gG)\rightarrow S(\gG)$ 
such that 
$
  f_{c}(X_1\cdot\cdot\cdot X_n) 
= c^n
  X_1\cdot\cdot\cdot X_n
$ 
for all $X_1,..., X_n\in\gG$. We have 
$
      f_{c}^{-1}
\circ \Phi^a_{c}
\circ f_{c}
    = \Phi^a_1
$ for any $a\in\gG$. \ \halmos

\section{A more general equation}
\label{sec22}

Let $\campo$ be a field of characteristic zero,  $t$ and $u$ be two 
commutating variables. We  classify the 
triples of formal series 
$(\varphi , \psi , \rho)$ such that $\varphi, \psi, \rho\in \campo[[t]]$ 
and
\begin{eqnarray}
\label{genL}
 \varphi(t)
 \frac{\psi(t+u)-\psi(u)}{t}
+\frac{\varphi(t+u)-\varphi(t)}{u}
 \psi(u)
= \rho(t+u)
.
\end{eqnarray}
This is motivated by theorem \ref{prop21} and it is clear that equation 
(\ref{precL}) is a  particular case of  equation (\ref{genL}). 
The classification is contained in theorems \ref{tm1} et \ref{tm2}. 

\begin{Oss}
Applying the limit $t\rightarrow 0$ we get 
\begin{eqnarray}
\label{dejaL}
  \varphi(0)\psi^{'}(u) 
+ \frac{\varphi(u) - \varphi(0)}{u}
  \psi(u)
= \rho(u)
.
\end{eqnarray}
\end{Oss}

\begin{Oss}
\label{oss232}
Applying limits $t\rightarrow 0$ and $u\rightarrow 0$ to equation 
(\ref{genL}) we get
$$  \varphi(0)\psi^{'}(u) 
- \psi(0)\varphi^{'}(u) 
= \frac{   \varphi(0)\psi(u)
         - \psi(0) \varphi(u)}
       {u}
.
$$
Thus there exists $a\in\campo$ 
such that   
\begin{equation}
\label{25}
  \varphi(0)\psi(u) 
- \psi(0)\varphi(u) 
= a u
.
\end{equation}
\end{Oss}

\noindent 
As the formal series $\rho$ is determined by (\ref{dejaL}),  it is natural 
to 
ask  
if equation (\ref{genL}) can be reduced to an equation for the
couple 
$(\varphi, \psi)$. To get such an  equation we introduce 
$   p(t), q(t)
\in  \frac{\campo}{t}+\campo[[t]]
$ 
such that 
$$
      \varphi(t) = tp(t)
, \ \ \psi(t) = tq(t)
.
$$ 

\begin{tm} 
\label{tmImpL}
The pair $(p(t), q(t))$ gives  a solution of 
(\ref{genL}) if and only if
\begin{equation}
\label{secondaL}
  q^{'}(u)
  \{ p(t+u)-p(t)\}
= p^{'}(t)
  \{ q(t+u) - q(u)\}
.
\end{equation}
\end{tm}
\noindent\dimostrazione\
Equation (\ref{genL}) becomes 
$$
  p(t)
  \{ (t+u)q(t+u) - uq(u) \}
+ q(u)
  \{ (t+u)p(t+u) - tp(t) \}
= \rho(t+u)
.
$$
We recall that a function $f(t,u)$ is a function of $t+u$ if and only if
 $\frac{\partial}{\partial t} f(t,u) - \frac{\partial}{\partial u} f(t,u) 
= 0$. We apply this fact to 
$
  f(t,u)
:=p(t)
  \{ (t+u)q(t+u) - uq(u) \}
+ q(u)
  \{ (t+u)p(t+u) - tp(t) \}
$ 
and we get equation (\ref{secondaL}). \halmos
\newline
\newline
Formula (\ref{secondaL}) is very elegant. However, we will use it only 
trough the following remark.
\begin{Oss}
\label{rem38}
If the pair $(p, q)$ is a solution without poles of
(\ref{secondaL}) then 
$
  q^{'}(u)
  (p(u) - p(0))
= 0
$. In particular  $p$ or $q$ is constant.

\end{Oss}
\begin{tm}
\label{tm1}
All triples of series  $(\varphi, \psi, \rho)$ 
verifying (\ref{genL}) and $\varphi(0)\psi(0) = 0$,  are given by  the 
following 
list 
\newline
i) $(\varphi(t), \psi(t), \rho(t)) = ( \varphi(t), ct, c\varphi(t) )$, \ \ \ \ 
$c\in\campo$, $\varphi\in\campo[[t]]$, 
\newline
ii) $(\varphi(t), \psi(t), \rho(t)) = (ct, \psi(t), c\psi(t) )$, \ \ \ 
$c\in\campo$, $\psi\in\campo[[t]]$.
\end{tm}
\noindent\dimostrazione\ 
It is sufficient to consider the case $\psi(0)=0$. From identity 
(\ref{25})  we get 
$
  \varphi(0)
  \psi(u) 
= a u
$, with $a\in\campo$. 

Let $\varphi(0) \neq 0$, we get $\psi(u) = cu$ with  $c\in\campo$. Equation 
(\ref{genL}) is verified with $c\varphi(t) = \rho(t)$, so we have a triple 
of type $i$.

Let  $\psi(0) = \varphi(0) = 0$. From remark \ref{rem38} we get that 
$q(u)$ or $p(u)$ is constant. It means that we get a triple of type $i$ or 
$ii$. \halmos
\newline

Now we treat the case $\varphi(0)\cdot \psi(0)\neq 0$.

\begin{Oss}
Let  $a,b\in\campo$. 
If $(\varphi(t), \psi(t), \rho(t))$  verifies the functional 
equation (\ref{genL}), then the triple 
 $(a\varphi(t), b\psi(t), a\cdot b \rho(t))$ verifies  (\ref{genL}).
\end{Oss}
\noindent 
Thus it is sufficient to look for series such that 
$\psi(0)=\varphi(0)=1$.

\begin{Oss}
Let  $a,b\in\campo$. If $(\varphi, \psi, \rho)$ verifies (\ref{genL}), 
then the  triple 
$
( \varphi(t)+at
, \psi + bt
, \rho(t) +a\cdot b t + a\psi(t) + b\varphi(t)
)
$ 
also verifies it.
\end{Oss}
\noindent
By this remark and by identity (\ref{25}), we can restrict  ourself to 
look for triples  
with $\varphi=\psi$,    $\varphi(0) = 1$, and $\varphi'(0)=0$. 

\begin{Oss}
\label{oss236}
Let  $a\in\campo$. 
If $(\varphi(t), \psi(t), \rho(t))$  verifies the functional 
equation (\ref{genL}),  then 
 $(\varphi(at), \psi(at), a \rho(at))$ verifies (\ref{genL}).
\end{Oss}
\noindent
For any $c\in\campo$ we introduce the notation 
$$
  \theta_c(t)
= \sqrt{c}t\coth(\sqrt{c}t)
= 1 + \frac{c}{3}t^2 - \frac{c^2}{45}t^4+\cdot\cdot\cdot
\in\campo[[t]]
.
$$
In particular $\theta_0(t)=1$.

\begin{lemma}
\label{lem231}
Let  $c\in\campo$. There exists exactly  one triple 
$(\psi_c, \psi_c, \rho_c)$ such that  
$\psi_c(t) = 1 + \frac{c}{3}t^2+ o(t^2)$. It is given by 
%Moreover  $\psi_c(t)=\psi_c(-t)$, 
$\psi_c(t) = \theta_{c}(t)$ and 
$\rho_c(t) =  ct$.
\end{lemma}
\noindent\dimostrazione\ 
We consider the left hand side of (\ref{genL}). As its derivative 
by $t$ must be equal to its  derivative by $u$ (see the proof of theorem 
\ref{tmImpL}),  
if  $\varphi= \psi$ we get 
\begin{eqnarray*}
  \psi'(t)
  \frac{\psi(t+u)-\psi(u)}{t}
+ \psi(t)
  \frac{\partial}{\partial t}
  \left( \frac{\psi(t+u)-\psi(u)}{t} \right)
+ \psi(u) 
  \frac{\psi'(t+u)-\psi'(t)}{u}
\\
= \psi(t)
  \frac{\psi'(t+u)-\psi'(u)}{t}
+ \psi'(u)
  \frac{\psi(t+u)-\psi(t)}{u}
+ \psi(u)
  \frac{\partial}{\partial u}
  \left( \frac{\psi(t+u)-\psi(t)}{u} \right)
.
\end{eqnarray*}
As $\psi(0)=1$ and $\psi'(0)=0$, the limit $t\rightarrow 0$ gives
\begin{eqnarray}
\label{l214}
- \frac{1}{2}
  \psi{''}(u)
= \frac{ \psi{'}(u) }{u}
  (\psi(u)-1)
- \psi(u)
  \frac{\psi(u)-1}{u^2}
.
\end{eqnarray}
Substituting $\psi(t) = 1+\frac{c}{3}t^2+\sum_{k=3}^\infty c_kt^k$  we 
get
$$
  c_{k+2}
= \frac{-2}{k(k+3)}
  \left(
     \frac{k c_{k}c}{3}
  +  \sum_{p=3}^{k-1}
     c_p c_{k-p+2}(k-p+1)
  \right)
, \ \ k\geq 1
.
$$
This formula gives $c_3= 0$. By induction,  all 
coefficients $c_{2j+1}$, $j\geq 2$ are zero.  The series $t\coth(t) = 1+ 
\frac{1}{3}t^2+\cdot\cdot\cdot$ is a solution of equation (\ref{l214}). By 
remark 
\ref{oss236} also  $ct\coth(ct)= 1+ \frac{c}{3}t+\cdot\cdot\cdot$ 
is a solution of equation (\ref{l214}). \halmos

%\begin{Oss}
%Let $c\in\campo$. By uniqueness,   $\psi_{c^2}(t) = \psi_1(ct)$ and 
%$\psi_{-c^2}(t) = \psi_1(ict)$. 
%\end{Oss}

\begin{tm}
\label{tm2}
The list of  solutions of (\ref{genL}) 
verifying  $\psi(0)\varphi(0) \neq 0$ is 
$$
\begin{array}{ll}
\left\{
\begin{array}{l}
  \varphi(t)
= a\theta_c(t) + bt 
\\
  \psi(t) 
= d\theta_c(t)
+ et
\\
  \rho(t) 
= (ae + bd)
  \theta_c(t)
+ (adc +be) t
\end{array}
\right.
\end{array}
$$
with $a,b,c,d,e\in\campo$ and $a\cdot d\neq 0$.
\end{tm}

%\begin{Oss}
%Let $\campo=\Razionali$. Triples of type $(\varphi, \varphi,-\varphi)$ 
%have $\varphi(t) = \varphi_{\frac{1}{2c}}(t)$ with $c\in\Razionali$  
%\end{Oss}

\begin{Oss}
\label{oss47}
We have 
\begin{equation}
\label{asL}
  \varphi_d
+ \varphi_{-d}
= \varphi_0
\equiv -t
, \ \ \forall d\in\campo^\times 
\end{equation}
and 
$
 \theta_{c}(t)
\equiv 
  \sqrt{c}
  \left(
  \varphi_{\frac{1}{2\sqrt{c}}}(t) 
- \varphi_{-\frac{1}{2\sqrt{c}}}(t) 
  \right)
= 2
  \sqrt{c}
  \left(
  \varphi_{\frac{1}{2\sqrt{c}}}(t) 
- \frac{1}{2}
  \varphi_0(t)
  \right)
$. 
\end{Oss}
%\noindent
%

\subsection{Application to a direct product}

Let $\campo$ be a commutative superring and $\gG$ a Lie 
$\campo$-superalgebra.

\begin{lemma}
\label{lemma71}
Let   $\campo\supseteq\Razionali$ and $g, h\in\campo^\times_0$.
The couple $(\Phi_g, \Phi_h)$ is composed of commuting   representations for 
any Lie $\campo$-superalgebra, if
and only if $g= -h$. 
\end{lemma}
\noindent\dimostrazione\ 
By theorem \ref{tm222}, $\Phi_g$ and $\Phi_h$ commute for any Lie 
$\campo$-superalgebra if and only if $(\varphi_g, \varphi_h, 0)$ verifies 
equation (\ref{genL}). As the formal series $\varphi_g$ and $\varphi_h$ 
have invertible coefficients, this is equivalent to the fact that 
$(\varphi_g, \varphi_h, 0)$  is one of the solutions given by theorem 
\ref{tm2}. Using that 
$$
  \varphi_c(t) 
= c
  \varphi_1\left( \frac{t}{c}
           \right)
, \ \ \forall c\in\campo_0^\times
$$
we see that $\varphi=\varphi_g$ if and only if $a=g$, $b=-\frac{1}{2}$, 
$c=\frac{1}{4g^2}$. In particular,  $(\varphi_g, \varphi_h, 0)$  is a 
solution of equation (\ref{genL}) if and only if 
$$
  0
= -\frac{1}{2}
  (g+h)\theta_c
+ \frac{1}{4}
  \left( \frac{h}{g}+1\right) t
.\  \halmos
$$

By denote by $\gG\times\gG$ the direct product of $\gG$ with its-self.
Let $\rho : \gG\times \gG \rightarrow
Hom(S(\gG), S(\gG))$ be a representation, it decomposes into the
sum of
two commuting representations $\rho_1, \rho_2 :\gG\rightarrow
Hom(S(\gG), S(\gG))$ such that
$\rho(a_1, a_2) = \rho_1(a_1) + \rho_2(a_2)$ for each
$(a_1,a_2)\in\gG\times \gG$. We write $\rho = (\rho_1, \rho_2)$. 
Theorem \ref{tm222} and lemma \ref{lemma71} give 

\begin{tm}
For any $\campo$-Lie 
superalgebra $\gG$, $ (\Phi_0, 0)$  and $(0, \Phi_0)$   are 
 representations by coderivations of $\gG\times \gG$ in 
$S(\gG)$.
\end{tm}

\begin{tm}
Assume that $\campo\supseteq\Razionali$. For any $\campo$-Lie 
superalgebra $\gG$, $ (\Phi_c, 0)$  and $(0, \Phi_c)$ whit 
$c\in\campo_0^\times\cup\{0\}$,  
$(\Phi_d, \Phi_{-d})$ with $d\in\campo_0^\times$, are 
 representations by coderivations of $\gG\times \gG$ in 
$S(\gG)$.
\end{tm}

\section{Lie algebras}
\label{sec23bis}

In this paragraph we consider the case of Lie algebras over a 
$\Razionali$-algebra. This means that   we assume that $\campo=\campo_0$
is a commutative 
ring, and $\gG$ a Lie $\campo$-algebra.

We have 
considered  coderivations associated to vector fields on $\gG$ of type 
\begin{eqnarray*}
\varphi^a =\varphi(\ad x)(a)
\end{eqnarray*} 
with $a\in\gG$, $\varphi\in\campo[[t]]$,  $x\in\gG_x$ the generic point of  
$\gG$. We have seen in remarks \ref{oss11} and \ref{rem219} that  
$\varphi^a$ and the corresponding coderivation $\Phi^a$
 satisfy a functorial property.  
In this paragraph we  prove the converse. 

We look for the family of morphisms of $\campo$-modules  
$
        F_\gG: S(\gG)\otimes\gG\rightarrow \gG
$ 
defined for all $\campo$-Lie algebra $\gG$ and such that, for all
morphisms of Lie $\campo$-algebras $f:\gG\rightarrow\hG$  
the diagram 
\begin{equation}
\label{diagramma}
\begin{array}{lll}
S(\gG)\otimes\gG
&\stackrel{F_\gG}{\longrightarrow} &\gG
\\
\downarrow^{\tilde f\otimes f}&& \downarrow^f
\\
S(\hG)\otimes\hG
&\stackrel{F_\gG }{\longrightarrow} 
&\hG
\end{array}
,
\end{equation}
where $\tilde f:S(\gG)\rightarrow S(\hG)$ is the  
algebra-morphism induced by $f$,  commutes.

\begin{tm}
\label{tm251}
For each $n\in\Naturali$, there exists 
$c_n\in\campo$ such that 
$$
  F_\gG(X_1\cdot\cdot\cdot X_n\otimes a)
= 
  ( c_n(\ad x)^n(a) )
  (X_1\cdot\cdot\cdot X_n)
$$
for any   $X_1,..., X_{n},a\in\gG$.
\end{tm}
\noindent\dimostrazione\ 
We consider  the free Lie $\campo$-algebra $\hG$ with 
generators 
$x_1,...x_{n+1}$. Let $Y:=F_\gG(x_1\cdot\cdot\cdot x_n\otimes 
x_{n+1}) $,  $t\in\Razionali$. We fix  $i\in\{1,...,n+1\}$. By the universal 
property of free Lie algebras, the map 
$$
        f_{t,i}: X_j 
\mapsto \left\{
        \begin{array}{ll}
        X_j,   &j\neq i
        \\
        X_it, &j=i
        \end{array}
\right.
, \ \ \forall j=1,...,n+1
$$
extends to a morphism of Lie $\campo$-algebras 
$\tilde f_{t,i}: \hG\rightarrow \hG$. 
As the  diagram (\ref{diagramma}) associated to this map commutes, 
we get $Yt = \tilde f_{t,i}(Y)$. 
We write $Y= \sum_{n} Y_{n,i}$ where $Y_{n,i}$ is a bracket containing $n$ 
times $x_i$, so $f_{t,i}(Y) \equiv \sum_n Y_{n,i} t^n$. To get 
$Yt = \sum_n Y_{n,i} t^n$, we need  
$\sum_{n\neq 1} Y_{n,i}(t - t^n)=0$ for any $t\in\Razionali$. 
As the family $\{Y_{n,i}|Y_{n,i}\neq 0, n\geq 0\}$ 
is free (see {\bf\cite{Bou2}}, prop. 10, page 26), 
we get 
that  
$Y=\sum_{i=1}^{n+1} Y_{1,i}$. This is true for any 
 $i\in\{i,..., n+1\}$, so $Y$ is a linear combination of brackets of $n$ 
elements, exactly elements $x_1,..., 
x_{n+1}$.  \newline
Using  the Jacobi identity and the fact that the bracket of a Lie algebra  
is antisymmetric, we show that $Y$ is a linear 
combination of 
$\ad x_{s(1)} \circ\cdot\cdot\cdot\circ\ad x_{s(n)}(x_{n+1})$, with 
$s\in\Sigma_{n}$. Let 
$ 
  Y 
= \sum_{s\in\Sigma_n}
  c_s
  \ad x_{s(1)} \circ\cdot\cdot\cdot\circ\ad x_{s(n)}(x_{n+1})
$, 
with $c_s\in\campo$. 
As a permutation  $u$ of $\{x_1,..., x_{n}\}$ extends to a morphism 
$g_u$ of Lie $\campo$-algebras, from the commutative diagrams 
(\ref{diagramma}) for  $g_u$ we get 
$ 
  \sum_{s\in\Sigma_n} 
  (c_s- c_{s\circ u})
  \ad x_{s(1)} \circ\cdot\cdot\cdot\circ\ad x_{s(n)}(x_{n+1})
=0
$. By the properties of free Lie algebras (see {\bf \cite{Bou2}}, prop. 10 
page 
26) we get that the family 
$
\{  \ad x_{s(1)} \circ\cdot\cdot\cdot\circ\ad x_{s(n)}(x_{n+1})
 |  s\in\Sigma_n
\}
$ is free. 
In particular  $c_s- c_{s\circ u}=0$ for any 
$s\in\Sigma_n$. As 
it is true for any permutation $s$, we get $c_s=c_{id}$ for any 
$s\in\Sigma_n$. We denote $c_{id}$ by $c_n$, so 
 $
  Y
= c_n
  \sum_{j\in\Sigma_n}
  \ad x_{j(1)} \circ\cdot\cdot\cdot\circ\ad x_{j(n)}(x_{n+1})
$.

Let $f$ be a map $\{x_1,..., x_n,  x_{n+1}\}\rightarrow \gG$ in a 
Lie $\campo$-algebra. From the universal property of free Lie 
algebras, $f$ extends 
to a morphism of Lie $\campo$-algebras still noted $f$. Let 
$a:=f(x_{n+1}), f(x_i)=:X_i$. The commutative diagram for $f$ gives 
\begin{eqnarray*}
&&F_\gG(X_1\cdot\cdot\cdot X_n\otimes a)
= c_n
  \sum_{j\in\Sigma_n}
  \ad X_{j(1)} \circ\cdot\cdot\cdot\circ\ad X_{j(n)}(a)
. \  \halmos
\end{eqnarray*}

\begin{Oss}
If $\campo=\campo_0$, the previous theorem is not valid for a Lie 
$\campo$-superalgebra $\gG$. For example, if $\theta_\gG=\theta$ is the 
map such 
that $\theta|_{\gG_0} = id $ and $\theta|_{\gG_1} = -id $ then 
$\gG\ni a\mapsto (\ad x)^n(\theta(a))$ has the functorial property 
expressed in diagram (\ref{diagramma}).
\end{Oss} 

As an application of theorem \ref{tm251} we get the following theorems 

\begin{tm}
\label{tm72}
Assume that $\campo=\campo_0\supseteq\Razionali$ is a field and 
$\gG=\gG_0$ is a Lie $\campo$-algebra. 
All universal  representations by coderivations 
$\gG\rightarrow Hom(S(\gG), S(\gG))$ are:  the zero representation, 
$\Phi_c$ with 
$c\in\campo$. 
\end{tm}

\begin{tm}
\label{tm63}
Assume that $\campo=\campo_0\supseteq\Razionali$ is a field and 
$\gG=\gG_0$ is a Lie $\campo$-algebra. 
All universal  representations by coderivations 
$\gG\times\gG\rightarrow Hom(S(\gG), S(\gG))$ are:  $(0, 0)$, $ (\Phi_c, 
0)$, $(0, 
\Phi_c)$,  $(\Phi_d, \Phi_{-d})$
with $c\in\campo$ and $d\in\campo\setminus\{0\}$.
\end{tm}

\section{The Poincar{\'e}-Birkhoff-Witt theorem}
\label{sec23}

Let $\campo$ be a commutative superring and  $\gG$ be a 
Lie $\campo$-superalgebra. We assume that 
 $\frac{1}{2}\in\campo$ or that 
$\campo=\campo_0$ and $\gG=\gG_0$.

We recall that the enveloping algebra $U(\gG)$ 
is defined as the quotient of the tensor algebra $T(\gG)$ by  the  ideal 
$J$ generated by 
$\{ a\otimes b 
  - (-1)^{p(a)p(b)}
    b\otimes a
  - [a,b]
    |  a,b\in\gG
\}$. 
The inclusion of $\gG$ in $T(\gG)$ gives a map $j:\gG\rightarrow U(\gG)$. 
Let $gr(U(\gG))$ be the graded module of $U(\gG)$ associated to the 
filtration 
$\{U_i\}_{i\geq 0}$ with $U_0=\campo$ and  $ U_i$ the $\campo$-module generated by 
$\{  j(X_1)\cdot\cdot\cdot j(X_l) 
   |  l\leq i, X_1,..., X_l\in\gG
\}$. The hypotheses give that it is a commutative superalgebra.

\begin{Oss}
If our assumptions  are not verified, $gr(U(\gG ))$ might not 
be commutative. For 
example let us consider $\gG=\Z e$ with odd $e$ and $[e,e]=0$. 
%$S^2(\gG)$ we have $e^2 = 0$, 
As $j(e)^2\notin U_1(\gG)$, $gr(U(\gG))$ is not commutative. 
\end{Oss}
\noindent
By the universal property of  symmetric algebras,  $j$ extends to  the 
algebra-morphism 
$$
\tilde j: S(\gG)\rightarrow gr(U(\gG))
$$
such that 
$X_1\cdot\cdot\cdot X_n\mapsto j(X_1)\cdot\cdot\cdot j(X_n) \mod U_{n-1}
$, for any  $n\in\Naturali$ and $X_1,..., X_n\in\gG$. This map is onto.

\begin{Def}
We say that $\gG$ verifies the {\it weak Poincar{\'e}-Birkhoff-Witt theorem} if 
$\tilde j$ is bijective. 
\end{Def}
\noindent
Before giving the next definition we recall that a  map 
$f: \gG\rightarrow \gG$ is said to be an {\it automorphism} if it 
is an invertible  morphism  of Lie $\campo$-superalgebras. 
Such a map induces two isomorphisms of $\campo$-superalgebras: 
$\tilde f :S(\gG)\rightarrow S(\gG)$ 
and 
$ \overline f: U(\gG)\rightarrow U(\gG)$.

\begin{Oss}
A derivation $g:\gG\rightarrow \gG$ extends to  derivations  
$ g_1: U(\gG)\rightarrow U(\gG)$ and 
$g_2: S(\gG)\rightarrow S(\gG)$. Moreover,  $g_1$ and $g_2$ are also 
two coderivations.
\end{Oss}

\begin{Def}
\label{def231}
We say that $\gG$ verifies the {\it strong Poincar{\'e}-Birkhoff-Witt
theorem} if it does 
exist an isomorphism $\rho\in Hom(S(\gG), U(\gG))$   such that
\newline
i) $\rho(S^n(\gG))\subseteq U_n(\gG)$ for any $n\in\Naturali$,
\newline
ii) the associated  graded map $gr(\rho)$ is  $\tilde j$,
\newline
iii) $\rho$ commutes with any derivation of $\gG$ and any  automorphism of 
$\gG$.
\end{Def}

\begin{Oss}
Let $n\geq 1$. From $i$ and $ii$ we have 
$ 
  \rho(S^n(\gG))\oplus U_{n-1}
= U_n
$. 
From ii we get that $\rho( S^n(\gG) )$ is stable by any derivation or 
automorphism of 
$\gG$. In particular $S(\gG)$ and $U(\gG)$ are isomorphic for the adjoint 
representation. 
\end{Oss}
\noindent
Now we suppose also that 

\begin{assume}
\label{hip253}
$\campo\supseteq\Razionali$ or  $\gG$ is 
$N$-nilpotent with $N\geq 2$ and
$\frac{1}{2},...,\frac{1}{N}\in\campo$. 
\label{as1}
\end{assume}
\noindent
From theorems \ref{tm222} and \ref{tm216} we have a 
representation $\Phi_1:\gG\rightarrow Hom(S(\gG), S(\gG) )$. By
the universal 
property of enveloping 
algebras, it extends to an algebra-morphism 
$
   \Phi : U(\gG)
\rightarrow Hom(S(\gG), S(\gG) )
$ such that $\Phi_1=\Phi\circ j$. From $\Phi$ we construct the map 
$
\sigma : U(\gG)
\rightarrow  S(\gG)
$, 
called the {\it symbol map} and defined by 
$$ 
   \sigma(u)
:= \Phi(u)(1)
, \ \forall u\in U(\gG)
.
$$ 

\begin{ex}
\label{ex71}
For any $a_1,a_2, a_3\in \gG$ we have 
$$
\begin{array}{l}
  \sigma(1) 
= 1
\\
  \sigma(j(a_1)) 
= a_1
\\
  \sigma(j(a_1) j(a_2)) 
= a_1\cdot a_2 + \frac{1}{2}[a_1, a_2]
\\
  \sigma( j(a_1) j(a_2) j(a_3) )
= a_1\cdot a_2\cdot a_3
+ \frac{1}{2}
  ( a_1\cdot [a_2, a_3]  
        + [a_1, a_2]\cdot a_3
        + (-1)^{p(a_1)p(a_2)}
          a_2\cdot [a_1, a_3]
  )
\\
\phantom{bbbbbbbbbbbbbbbbbbbbbb}
+ \frac{1}{12}
  \left\{ 
        -(-1)^{p(a_2)p(a_1)}
          [a_2, [a_1, a_3]] 
        + [[a_1, a_2], a_3]
  \right\}
+ \frac{1}{4}
  [a_1, [a_2, a_3]] 
.
\end{array} 
$$
\end{ex}

\begin{lemma}
\label{lemma52}
Let $\campo$ be any commutative superring and $\gG$ any
Lie $\campo$-superalgebra. If   $\lambda\in\campo[[z]]$ with 
$\lambda(0)=1$,  the 
coderivations  corresponding to $\lambda$  have the property
$$
\Lambda^{a_n}\circ\cdot\cdot\cdot\circ\Lambda^{a_1}
  ( 1)
- a_n\cdot\cdot\cdot a_1
\in \oplus_{j=0}^{n-1}S^j(\gG)
, \ \ \forall a_1,..., a_n\in\gG
.
$$
\end{lemma}
\noindent\dimostrazione\ 
If $n=1$ the theorem is evident. As
$
  \Lambda^{a_n}\circ\cdot\cdot\cdot\circ\Lambda^{a_1}
  ( 1)
= \Lambda^{a_n}
  \left( \Lambda^{a_{n-1}}\circ\cdot\cdot\cdot\circ\Lambda^{a_1}(1)
  \right)
$,  
by induction there exists $p_n\in \oplus_{j=0}^{n-2}S^j(\gG)$ such 
that 
$ \Lambda^{a_{n-1}}\circ\cdot\cdot\cdot\circ\Lambda^{a_1}(1) 
= a_{n-1}\cdot\cdot\cdot a_1 + p_n
$. 
As for any  $p\geq 0$
$$
        \Delta(S^p(\gG))
\subseteq \bigoplus_{j=0}^p
        S^j(\gG)\otimes S^{p-j}(\gG), 
$$
we have 
$         \Lambda^{a_n}(p_n)
\in \bigoplus_{j=1}^{n-1}
           S^j(\gG) 
$,  
it is sufficient to show that 
$
    \Lambda^{a_n}(a_{n-1}\cdot\cdot\cdot a_1)
  - a_n\cdot\cdot\cdot a_1
\in \oplus_{j=0}^{n-1}S^j(\gG)
$.  
This identity follows using 
$
    \Delta(a_{n-1}\cdot\cdot\cdot a_1)
  - a_{n-1}\cdot\cdot\cdot a_1\otimes 1
\in \oplus_{j=0}^{n-2}
    S^j(\gG)\otimes S^{n-1-j}(\gG)
$, 
 the definition of $\Lambda$, the identity (\ref{phi1}). \halmos

\begin{tm}\footnote{See the historical note below}
\label{tm22}
%\label{tm42}
%
Assume hypothesis \ref{hip253}. Then 
\newline
i) the map $\sigma$ is invertible, 
\newline
ii)  $\gG$ verifies the strong 
Poincar{\'e}-Birkhoff-Witt theorem with $\rho=\sigma^{-1}$.
\end{tm}
\noindent\dimostrazione\ 
$i)$ 
By lemma  \ref{lemma52} the graded map 
 $gr(\sigma) : gr(U(\gG))\rightarrow S(\gG)$ 
is well-defined and onto: for any $n\in\Naturali$ we have 
$$
  gr(\sigma)
  ( a_1\cdot\cdot\cdot a_n + U_{n-1})
= \sigma(a_1)\cdot\cdot\cdot\sigma(a_n)
, \ \ \forall a_1,..., a_n\in\gG
. 
$$
The inverse of $gr(\sigma)|_{U_n/U_{n-1}}$ is $\tilde j|_{S^n(\gG)}$, so 
$gr(\sigma)$ is one-to-one. 
\newline
$ii)$ In $i$ we have seen, in particular, that $\tilde j =
gr(\sigma^{-1})$.
Let $f: \gG\rightarrow \gG$ be an automorphism of $\gG$. The fact that 
$
  \sigma\circ \overline f
= \tilde f\circ \sigma
$ follows from  remark \ref{rem219}.

Let $g$ be a derivation of $\gG$ and $a_1,..., a_n\in\gG$. We want to
show that 
$     g_2
\circ \sigma 
    = \sigma
\circ g_1
$, 
which means 
$$ 
\left(
      g_2
\circ \Phi^{a_1}
\circ \cdot\cdot\cdot
\circ \Phi^{a_n}
\right)
      (1)
    = \sum_{j=1}^n
      (-1)^{p(g)p(a_1+\cdot\cdot\cdot +a_{j-1})}
\left(
      \Phi^{a_1}
\circ \Phi^{g(a_j)}
\circ \Phi^{a_n}
\right)
      (1)
$$
for any $a_1,..., a_n\in\gG$. 
By induction, it is sufficient to show that 
$ [g_2, \Phi_1^{a}]
= \Phi^{g(a)}
$ for any $a\in\gG$.
By definitions  
$ g_2\circ \varphi_1^a
= \varphi_1^a\circ g_2
+ \varphi^{g(a)}
$. This gives 
$
      [g_2, \Phi_1^{a}]
    = 1\otimes \varphi_1^{a}
\circ ( g_2 \otimes 1 + 1\otimes g_2
      - \Delta\circ g_2
      )
    + \Phi^{g(a)}
$. 
The fact that $ g_2$ is a coderivation ends the proof. 
\halmos 
\newline
\newline
Let $\beta:= \sigma^{-1}$.

\begin{Oss} 
\label{oss244}
{\bf (Functorial property)}
\newline
Let $f:\gG\rightarrow \hG$ be a  morphism of 
Lie $\campo$-superalgebras.  By remark 
\ref{rem219} we get  a commuting diagram
$$
\begin{array}{ccc}
S(\gG)&\stackrel{\beta}{\longrightarrow} &U(\gG)
\\
\downarrow_{\tilde f}
&&\downarrow_{\overline f}
\\
S(\hG)&\stackrel{\beta}{\longrightarrow} &U(\hG)
\end{array}
.
$$ 
\end{Oss}
\noindent
To get formulas for $\beta$ we use the following lemma.

\begin{lemma}
\label{lemma41}
For any $n\in\Naturali$ and  $a\in\gG_0$
we have 
$
  (\Phi_1^a)^n(1)
= a^n 
$.
\end{lemma}
\noindent\dimostrazione\ 
If $n=1$ the statement is obvious. By induction 
$$
       (\Phi^a_1)^{n+1}(1)
\equiv \Phi^a_1\circ (\Phi^a_1)^n(1)
     = \Phi^a_1(a^n)
     = \sum_{j=0}^n
       {n\choose j}
       a^j\cdot\varphi^a(a^{n-j})
.
$$
From identity (\ref{labiii}) we get $\varphi^a(a^j)= 0$ for $j\geq 1$, so 
$ 
  (\Phi^a_1)^{n+1}(1)
= a^n\cdot\varphi^a(1)
= a^{n+1} 
$.
\ \halmos

\begin{Cor}
Let $n\in\Naturali$. 
\newline 
i) For any $a\in\gG_0$ we have 
$
  \beta(a^n)
= \beta(a)^n
= j(a)^n
$.
\newline
ii) For any $a_1,..., a_n\in\gG$, 
$
  n!
  \beta(a_1\cdot\cdot\cdot a_n)
= \sum_{s\in\Sigma_n}
  \alpha( a_{s(1)}, ..., a_{s(n)})
  \beta(a_{s(1)})\cdot\cdot\cdot\beta(a_{s(n)})
.
$
\end{Cor}
\noindent
From now on $\beta$ will be called the {\it symmetrization map}. If 
$\campo$  contains $\Razionali$, $\beta$ is the usual symmetrization 
map. If 
$\campo$ does not contain $\Razionali$, the 
previous corollary does not give an explicit formula for the 
symmetrization map. However, in principle we can compute 
$\beta(a_1\cdot\cdot\cdot 
a_n)$ 
(as in the following example) but we do not know a nice formula.

\begin{ex} 
Let $\campo=\Z/3\Z$, it contains $\frac{1}{2}$. Let 
$\gG$ be a $2$-nilpotent Lie superalgebra over $\campo$. For 
each $a\in\gG$, $\Phi^a = a^L + \frac{1}{2}\ad a$.  
Let $a_1,a_2, a_3\in\gG$. From example \ref{ex71} we get 
\begin{eqnarray*}
&&\sigma(j(a_1)) = a_1
\\
&&\sigma(j(a_1)j(a_2)) = a_1\cdot a_2 + \frac{1}{2}[a_1, a_2] 
\\
&&  \sigma(j(a_1)j(a_2)j(a_3)) 
= a_1\cdot a_2\cdot a_3 
+ \frac{1}{2}
  \left\{
          a_1\cdot [a_2, a_3]
        + [a_1, a_2]\cdot a_3 
\right\}+
\\
&&
\phantom{bbbbbbbbbbbbbbbbbbbbbbbbbbbbbbbbbbbbb}
+ \frac{1}{2}
  (-1)^{p(a_2)p(a_1)}
  a_2\cdot [a_1,a_3]
.
\end{eqnarray*}
Applying $\beta$ we get 
\begin{eqnarray*}
&&    \beta(a_1) = j(a_1)
\\
&&    \beta(a_1\cdot a_2)
= j(a_1)j(a_2) - \frac{1}{2}\beta([a_1, a_2]) 
  = j(a_1)j(a_2) - \frac{1}{2}j([a_1, a_2] )
\\
&&   \beta(a_1\cdot a_2\cdot a_3) =
\\
&=&j(a_1)j(a_2)j(a_3 )
- \frac{1}{2}
  \beta
  \left(
          a_1\cdot [a_2, a_3]
        + [a_1, a_2]\cdot a_3 
        +(-1)^{p(a_2)p(a_1)}
         a_2\cdot [a_1,a_3]
\right)=
\\
&=& j(a_1)j(a_2)j(a_3) 
- \frac{1}{2}
  (j(a_1)j([a_2, a_3])
+ j([a_1, a_2])j(a_3)
  )+ 
\\
&&\phantom{bbbbbbbbbbbbbbbbbbbbb}
- \frac{1}{2}
  (-1)^{p(a_2)p(a_1)}
  j(a_2)j([a_1,a_3])
.
\end{eqnarray*}
\end{ex}

\begin{Oss} 
{\bf (Historical note)}
\newline
In the literature you can find proofs of the fact that the 
symmetrization $\beta$ 
is an isomorphism of $\campo$-modules for  $\campo=\campo_0\supseteq 
\Razionali$, $\gG$ a Lie $\campo$-algebra ({\bf \cite{Coh}},  
{\bf \cite{Bou2}} exercise 16, page 78) or a Lie superalgebra ({appendix 
of \bf \cite{Qui}}). 
Even in these particular cases, our proof is different and more direct. 
In particular, we do not have to consider first the special case of free 
Lie superalgebras.
 
The case of $N$-nilpotent Lie superalgebras was  known only for $N=2$. 
It 
was 
proved by M. El-Agawany and A. Micali  (see {\bf \cite{ElM}}).

Before {\bf\cite{Coh}} the strong theorem was  
known for some class of Lie algebras. For example, if $\campo$ is a field 
of characteristic zero, it is due to Poincar{\'e} (see {\bf\cite{Tut}}). 
Examples 
of Lie algebras not verifying the weak Poincar{\'e}-Birkhoff-Witt theorem  
are given in {\bf \cite{Sir}}, {\bf 
\cite{Car}}, {\bf \cite{Coh}}. 

If $\campo$ is a field (or more generally if $\gG$ is a free 
$\campo$-module), then 
any Lie $\campo$-algebra verifies the weak Poincar\'e-Birkhoff-Witt 
theorem (see for example {\bf\cite{Bou2}}). 

This is also true for a Lie $\campo$-superalgebra if $2$ is invertible in 
$\campo$ 
(see {\bf \cite{BMP}}). However, if $\campo$ is a field of finite 
characteristic, the strong Poincar\'e-Birkhoff-Witt theorem is usually not 
satisfied. 

\end{Oss}

\subsection{Universal representations in the enveloping algebra}

We still assume hypotheses \ref{as1}. By theorem \ref{tm22} we can 
transport each 
coderivations  $\Phi_c$, $c\in\campo^\times_0\cup\{0\}$, on $U(\gG)$.

We recall that $U(\gG)$ is equipped of a natural 
comultiplication $\Delta^{'}$, such that for $a\in\gG$ we have 
$\Delta^{'}(j(a))= 1\otimes j(a) + j(a)\otimes 1$.

\begin{tm}
\label{OSS1}
The symmetrization map verifies 
$     \Delta^{'}
\circ \beta
    = (\beta\otimes \beta)
\circ \Delta
$.  
In particular for any  $a\in\gG$ and $c\in\campo_0^\times\cup\{0\}$, 
$\beta\circ \Phi_c^a\circ \beta^{-1}$ is a coderivation of $U(\gG)$.
\end{tm}
\noindent\dimostrazione\
We consider the map $\gG\rightarrow \gG\times\gG$ such that 
$X\mapsto (X,X)$. It induces the comultiplication 
over
$S(\gG)$ and $U(\gG)$, so remark (\ref{oss244}) ends the proof. \halmos
\newline
\newline
Let $a\in\gG$. In $U(\gG)$ we 
have $\ad j(a) = j(a)^L - j(a)^R$ (see notation \ref{not21}). 

\begin{tm}
\label{MenomaleL}
Let $a\in\gG$, 
\newline
i)  $\beta^{-1}\circ \ad j(a)\circ \beta = \Phi^a_0$
\newline
ii) $\beta^{-1}\circ j(a)^L\circ \beta = \Phi^a_1$
\newline
iii) $\beta^{-1}\circ j(a)^R\circ \beta = -\Phi^a_{-1}$.
\end{tm}
\noindent\dimostrazione\ 
$i)$
The map  $\gG\ni X\mapsto [a,X]$ is a derivation of $\gG$, it 
extends to the derivations $\ad j(a)$ and $\Phi_0^a$. Theorem \ref{tm22} 
gives the identity $i$. 
\newline
$ii)$ Let  $W\in S(\gG)$. To show that 
$
\sigma\left( j(a)\cdot \beta(W) \right)
= \Phi^a_1(W)
$
we only need to recall that by definitions we have 
$
       \sigma\left( j(a)\cdot \beta(W) \right)
     = \Phi_1^a\circ \sigma\left( \beta(W)\right)
\equiv \Phi_1^a(W)
$.
\newline
$iii)$ As $\ad j(a) = j(a)^L - j(a)^R$ in $U(\gG)$, cases $i$ and $ii$ give 
$
  \beta^{-1}\circ j(a)^R\circ \beta 
= \Phi^a_1
- \Phi^a_0
$. 
By identity (\ref{asL}), the coderivation $\Phi^a_1- \Phi^a_0$ is equal to
$-\Phi^a_{-1}$.\ \halmos

\begin{Oss}
The map $a\mapsto \beta\circ\Phi_c^a\circ\beta^{-1}$ interpolates the 
regular left 
representation $a\mapsto j(a)^L$ ($c=1$) and the regular right  
representation $a\mapsto -j(a)^R$ ($c=-1$).
\end{Oss}

\begin{tm}
Let $\campo=\campo_0$ be a field  of characteristic zero and $\gG=\gG_0$ 
a Lie $\campo$-algebra. 
All 
universal representations $\gG\rightarrow Hom( U(\gG), U(\gG) )$ by 
coderivations  are equivalent to the zero representation, or to the 
adjoint representation, or to the regular left representation.
\end{tm}
\noindent\dimostrazione\ 
Let $F:\gG\rightarrow Hom(U(\gG), U(\gG))$ be a representation by 
coderivations. We assume that $F$ is not the zero representation. 
By theorem \ref{OSS1}, for any $a\in\gG$, $G(a):=\beta^{-1}\circ 
F(a)\circ\beta$ 
is a coderivation of $S((\gG)$. In particular 
$G$ is a representation by 
coderivations of $\gG$ in $S(\gG)$. Using that $\campo$ is a field and 
using theorem \ref{tm251},  we get that $G$
is one of the representations given in     theorems \ref{tm213} and 
\ref{tm222}. From 
theorem \ref{cor213} we 
get that $G$ is equivalent to 
$\Phi_1$ or $\Phi_0$. By theorem \ref{MenomaleL}, $\gG\ni a\mapsto 
G(a)$ 
is equivalent to  $\gG\ni a\mapsto \ad j(a)$  or $\gG\ni a\mapsto j(a)^L$. 
\halmos
\newline
\newline
Using  theorems \ref{cor213}, \ref{tm63},  and \ref{MenomaleL},  we get

\begin{tm}
Let $\campo=\campo_0\supseteq\Razionali$ be a field and $\gG=\gG_0$ a 
$\campo$-Lie algebra. We have 5 classes of 
equivalence for non-zero  universal representations by
coderivations  $\gG\times\gG\rightarrow Hom(U(\gG), U(\gG))$:
\begin{eqnarray*}
&&  \gG\times\gG\ni(a,b)\mapsto \alpha\ad a+ (1-\alpha)\ad b
, \ \ \ \ \ \ \ \alpha\in\{0,1\}
\\
&&\gG\times\gG\ni (a,b)\mapsto \alpha j(a)^L-(1-\alpha) j(b)^R
, \ \ \ \ \ \ \ \alpha\in\{0,1\}
\\
&&\gG\times\gG\ni (a,b)\mapsto j(a)^L - j(b)^R
.
\end{eqnarray*}
\end{tm}

\begin{Oss}
Let $\campo$ a field of  characteristic zero. We denote by  
$\pi_1: S(\gG)\rightarrow \gG$ the projection over $\gG$ and we put 
$P = \pi_1\circ\beta^{-1}$. By theorem  \ref{lemma11}, 
the parts  $ii$ and $iii$ of theorem 
\ref{MenomaleL} are equivalent to 
$
P\circ a^L = \varphi_1^a\circ\beta^{-1}
, \ \
P\circ a^R = -\varphi_{-1}^a\circ\beta^{-1}
$. 
In  {\bf\cite{Sol}} and {\bf \cite{Hel}} we can find a formula for  $P$. 
\end{Oss}

\begin{Oss}
In   {\bf \cite{Ber}} and {\bf \cite{Ras}} we can find the   formula 
for $\beta^{-1}\circ a^L\circ\beta$ which is in theorem 
\ref{MenomaleL}.
\end{Oss}

\end{document}